\documentclass[12pt]{article}
\usepackage[utf8]{inputenc}
\usepackage[left=18mm, right=18mm, top=20mm, bottom=20mm]{geometry}
\usepackage{amsmath,amsfonts,amssymb,amsthm}
\usepackage{hyperref}
\usepackage{graphicx}
\usepackage{epstopdf}
\usepackage{color}
\usepackage{float}
\usepackage{hyperref}
\usepackage{endnotes}
\usepackage{ragged2e}
\theoremstyle{definition}
\usepackage{amssymb}
\usepackage{amsmath}
\usepackage{amsthm}
\usepackage{float}
\usepackage{booktabs}
\usepackage{graphicx}
\usepackage{epstopdf}
\usepackage{lineno}

\newcounter{bla}

  \theoremstyle{definition}

\usepackage{times}
\numberwithin{equation}{section}
\newtheorem{thm}{Theorem}[section]
\newtheorem{lem}{Lemma}[section]

\providecommand{\keywords}[1]
{
	\small	
	\textbf{\textit{Keywords---}} #1
}
\providecommand{\MSC}[1]
{
	\small	
	\textbf{\textit{MSC Subject Classification---}} #1
}
\title{Simple smoothness indicator WENO-Z scheme for hyperbolic conservation laws}
\author{Samala Rathan$^{1,2}$\thanks{\textit{Email:rathans.math@iipe.ac.in}} , G. Naga Raju$^{2}$\thanks{\textit{Email:gnagaraju@mth.vnit.ac.in}}, Ashlesha A. Bhise$^{2}$\thanks{\textit{Email:ashleshaabhise@gmail.com}}  \\
	\small $^{1}$\textit{Faculty of Mathematics, Indian Institute of Petroleum \& Energy-Visakhapatnam, India-530003} \\
	\small $^{2}$\textit{Department of Mathematics,
		Visvesvaraya National Institute of Technology, Nagpur, India-440010} \\
}
\date{} 
\begin{document}
	\maketitle	
	\begin{abstract}
	The advantage of WENO-JS5 scheme [ J. Comput. Phys. 1996] over the WENO-LOC scheme [J. Comput. Phys.1994] is that the WENO-LOC nonlinear weights do not achieve the desired order of convergence in smooth monotone regions and at critical points.  In this article, this drawback is achieved with the WENO-LOC smoothness indicators by constructing a WENO-Z type  nonlinear weights which contains a novel global smoothness indicator.  This novel smoothness indicator measures the derivatives of the reconstructed flux in a global stencil, as a result, the proposed numerical scheme could decrease the dissipation near the discontinuous regions. The theoretical and numerical experiments to achieve the required order of convergence in smooth monotone regions, at critical points, the essentially non-oscillatory (ENO), the analysis of parameters involved in the nonlinear weights like $\epsilon$ and $p$  are studied.  From this study, we conclude that the imposition of certain conditions on $\epsilon$ and $p$,  the proposed scheme achieves the global order of accuracy in the presence of an arbitrary number of critical points.   Numerical tests for scalar, one and two-dimensional system of Euler equations are presented to show the effective performance of the proposed numerical scheme.
	\end{abstract}
	\\
	\MSC{65M20, 65N06,  41A10.}
\section{Introduction}
The study of hyperbolic conservation laws 
\begin{equation}\label{e1}
\begin{aligned}
\frac{\partial \mathbf{u}}{\partial t}+ \frac{\partial f(\mathbf{u})}{\partial x} &= 0, \mathbf{x}\in \mathbb{R}^{d},(d \geq 1), t>0,\\
\mathbf{u}(x,0)&= \mathbf{u}_0(x), \mathbf{x}\in \mathbb{R}^{d},
\end{aligned}
\end{equation} 
is one of the important topics  in the areas of gas dynamics, shallow water flows and magneto-hydro-dynamics(MHD). For the equation (\ref{e1}),  $\mathbf{u}=(u_{1},u_{2},......,u_{m})^{T}$ represents a conserved quantity which is a $m$-dimensional vector  and flux $f(\mathbf{u})$, is a vector-valued function of $m$ components, $x$ and $t$ denote the space and time variables respectively. 
It is well known that the analytical solutions are available only for a few model problems and thus,  numerical techniques play a important role in solving problems of practical interest. The vital remark in the solutions of hyperbolic conservation laws is  that even if the smooth initial data may give rise to discontinuities as the time is propagating. For resolving this scenario and to obtain a valid solution, many numerical techniques such as finite difference, finite volume and finite element techniques have been developed.
\newline
Among them, the essentially non-oscillatory(ENO) schemes \cite{HO97, HO87, osherSHU,Shu-osher1} and the weighted essentially non-oscillatory (WENO) schemes \cite{XD Liu, Jiang and shu7} are quite popular. As our interest is on WENO schemes, we briefly mention the details about these schemes. The WENO schemes first developed in $1994$ by Liu, Osher, and Chan \cite{XD Liu} in a finite-volume framework where the authors came up with an ingenious idea as such: instead of choosing the smoothest candidate stencil, a nonlinear convex combination of all the sub stencils is used which results overall, a high-order accurate scheme when it is compared to ENO schemes. The major contributions of this technique are the construction of the nonlinear weights and the smoothness indicators based on undivided differences. Later in $1996$ \cite{Jiang and shu7}, a finite difference WENO schemes are developed with the construction of new smoothness indicators, commonly known as WENO-JS (JS stands for Jiang \& Shu) schemes. Hereafter, we refer the finite difference WENO formulation with the smoothness indicators of \cite{XD Liu} as WENO-LOC scheme. The smoothness indicators of the WENO-JS scheme are the square sum of all the derivatives of $m$ local interpolating polynomials, the process leads to obtaining $(2m-1)^{th}-$order accuracy of the scheme in smooth regions.  These schemes are extended by Balsara and Shu in  \cite{balsara and shu8} to a WENO family up to  $11^{th}-$order accuracy. Besides, Gerolymous et al. \cite{Gerolymus} introduced a WENO family up to  $17^{th}-$order.  Balsara et al. \cite{Balsara1}  analyzed the WENO scheme presented in \cite{balsara and shu8} in a basis set formed by Legendre polynomials up to $9^{th}-$order which affords an equivalent formulation for the numerical fluxes, as a result, the smoothness indicators are in the compact form. And further, the smoothness indicators have been written as the sum of perfect squares which makes the method more efficient and also more accurate for certain benchmark problems. This procedure further  carried out in \cite{Balsara2} up to $17^{th}-$order. Henrick et al. \cite{henrick aslam powers5} studied the WENO-JS scheme and discovered that the WENO-JS nonlinear weights failed to recover the optimal order of accuracy at the critical points where the first-order derivative vanish but not the third-order derivative and observed that the scheme is sensitive with respect to the choice of $\epsilon$, the parameter used in the evaluation of nonlinear weights. To dissolve this issue and to achieve the required order of accuracy in presence of critical points, the authors altered the nonlinear weights through the construction of a mapping function which approximates the WENO convex combination intently to the optimal weights except at highly non-smooth regions. Another approach was adapted by Borges et al. \cite{Borges caramona10} where the author's designed global smoothness measurements for the fifth-order WENO scheme, dubbed as WENO-Z, which has the same accuracy as that of mapped WENO with the lower computational cost. Castro et al. \cite{Castroetal} extended WENO-Z schemes to the higher-order, which have computationally cheaper nonlinear weights than mapped WENO  through the construction of high-order smoothness indicators that can be obtained from the inexpensive linear combination of existing lower order smoothness indicators.  Many modified and improved versions of the WENO schemes can be seen \cite{Ha13, KIM et al, rathan, rathan1, serna, Fan14, Fan15, Feng, Wang, YC, BD, DP, RP}.
\par
 It is well known that the WENO schemes are quite popular from last two decades to approximate the solutions of the hyperbolic conservation laws through the smoothness  indicators developed in \cite{Jiang and shu7}, the authors,  Jiang and  Shu, modified  the smoothness indicators developed in   \cite{XD Liu}   as in smooth monotone regions and at the critical points the WENO-LOC scheme does not achieve the desired order of accuracy.  To resolve this, in this paper, we have constructed a WENO-Z type nonlinear weights with WENO-LOC smoothness indicators. A novel global smoothness indicator is devised by measuring the derivatives of the reconstructed flux through undivided differences, as a result, the numerical scheme could decrease the dissipation around the discontinuities. Further, the proposed numerical scheme achieves the sufficient condition and ENO property to gain the required order of accuracy in smooth regions and at critical points.  Several benchmark problems in the scalar, the system of one- and two-dimensional  Euler equations are performed to show the effective performance of the proposed numerical scheme. It is shown that the proposed WENO scheme provides improved behavior to the fifth-order WENO-LOC and fifth-order WENO-JS (WENO-JS5) schemes. Furthermore, the consistency analysis of the numerical scheme is developed and shown that the imposition of certain conditions on the weight parameters leads to achieve the desired global order of accuracy in the presence of the arbitrary number of critical points. 
\par
The rest of the paper is organized as follows. The detailed formulation of the WENO scheme with the WENO-LOC and WENO-JS5 schemes are given in Section 2.  In Section 3, the design of new nonlinear weights is proposed and performed the ENO property, accuracy test in smooth regions, near discontinuities and at critical points. Numerics have been performed for some benchmark problems like a scalar, one and two-dimensional Euler equations in Section 4. Concluding remarks are given in Section 5.
\section{Numerical Scheme}
In this section, for completeness we report the flux version of  fifth-order WENO schemes presented in \cite{Jiang and shu7} for hyperbolic conservation laws (\ref{e1}).
\subsection{WENO schemes}\label{sec1}
\par
Let $\{I_{i}\}_{i}$ with $I_i=[x_{i-\frac{1}{2}},x_{i+\frac{1}{2}})$ be the partition of computational domain  in space and let $x_i=\frac{1}{2}(x_{i+\frac{1}{2}}+x_{i-\frac{1}{2}})$ denote the center of the cell $I_i$ with the uniform cell length  $\Delta x=x_{i+\frac{1}{2}}-x_{i-\frac{1}{2}}$. The function value $f$ at the node $x_i$ is given by $f_i:=f(x_i).$ Moreover, we use the notation $u_i^n$ for the approximation to  $u$ at the grid point $(x_i,t^n)$ and   $t^n=n\Delta t$. For simplicity, we restrict our discussion to one-dimensional scalar  formulation of  (\ref{e1}),
\begin{equation}\label{fd}
u_t=-f(u)_x,
\end{equation}
 and the associated semi-discretized formulation  is
\begin{equation}\label{fhat}
\frac{d{u_{i}(t)}}{dt}=-\frac{1}{\Delta x}\left({\hat{f}_{i+\frac{1}{2}}-\hat{f}_{i-\frac{1}{2}}}\right)=:L(u),
\end{equation}
where $u_i(t)$ is the numerical approximation to the point value $u(x_i,.)$ and the numerical flux $\hat{f}$ is a function of $(r+s)$ arguments i.e., $\hat{f}_{i+\frac{1}{2}}=\hat{f}(u_{i-r},...,u_{i+s}).$  The  system of ODE's  (\ref{fhat}) can be obtained  by using the strong-stablity preserving Runge-Kutta methods \cite{Gottilieb}.
\par
The numerical flux function $\hat{f}$ in (\ref{fhat}) should be consistent with the physical flux $f,$ that is, $\hat{f}(u,...,u)=f(u)$ and should satisfy the Lipschitz continuity in each of its arguments, as a requirement for the applicability of Lax-Wendroff theorem \cite{RJL}.
\par
To compute the numerical flux $\hat{f}_{i \pm \frac{1}{2}},$ a function $h$ is defined implicitly (see Lemma 2.1 of \cite{Shu-osher1})
\begin{equation}\label{e4}
f(x):=f(u(x,.))=\frac{1}{\Delta x}\int_{x-\frac{\Delta x}{2}}^{x+\frac{\Delta x}{2}}h(\xi)d\xi.
\end{equation}
The differentiation of equation \eqref{e4} and evaluation at the point $x=x_i$ yields
\begin{equation}\label{e5}
\frac{\partial f}{\partial x}\bigg|_{x=x_i}=\frac{1}{\Delta x}\left(h_{i+\frac{1}{2}}-h_{i-\frac{1}{2}}\right),
\end{equation}
which indicates that the numerical flux $\hat{f}$  approximates $h$ at cell boundaries $x_{i\pm \frac{1}{2}}$ with high-order of accuracy, that is,
\begin{equation*}
\hat{f}_{i \pm \frac{1}{2}}=h(x_{i \pm \frac{1}{2}})+O\big(\Delta x^k\big),
\end{equation*}
where $k$ depending on the degree of interpolation.  The basic  observation  reveals that the  spatial derivative defined  in \eqref{fd} is exactly approximated by a conservative finite difference formula \eqref{e5} at the cell boundaries. Using equation \eqref{e5} in equation \eqref{fhat}, we have
\begin{equation}\label{e9}
\frac{du_i(t)}{dt}=-\frac{1}{\Delta x}\left(h_{i+\frac{1}{2}}-h_{i-\frac{1}{2}}\right)\approx -\frac{1}{\Delta x}\left(\hat{f}_{i+\frac{1}{2}}-\hat{f}_{i-\frac{1}{2}}\right).
\end{equation}
In order to ensure the numerical stability,  the flux $f(u)$ is splitted into two parts $f^{+}$ and $f^{-}$ such that
\begin{equation}\label{fs}
f(u)=f^{+}(u)+f^{-}(u),
\end{equation}
where $\frac{df^{+}(u)}{du}\geq0$ and $\frac{df^{-}(u)}{du}\leq0.$ Among many flux splitting methods, we use global Lax-Friedrichs splitting 
\begin{equation}
f^{\pm}(u)=\frac{1}{2}(f(u)\pm \alpha u),
\end{equation}
where $\alpha =\displaystyle\max_{u}|f{'}(u)|$ for its simplicity and capability to produce very smooth fluxes.  Let $\hat{f}^{+}_{i+\frac{1}{2}}$ and $\hat{f}^{-}_{i+\frac{1}{2}}$ be the numerical fluxes obtained from the positive and negative parts of $f(u)$ respectively and from (\ref{fs}), we have
\begin{equation}
\hat{f}_{i+\frac{1}{2}}=\hat{f}_{i+\frac{1}{2}}^{+}+\hat{f}_{i+\frac{1}{2}}^{-}.
\end{equation}
Now we describe only how $\hat{f}_{i+\frac{1}{2}}^{+}$ can be approximated
since $\hat{f}_{i+\frac{1}{2}}^{-}$ is symmetric to the positive part
with respect to $x_{i+\frac{1}{2}}.$ In the formulation
of $\hat{f}_{i+\frac{1}{2}}^{+}$,  for simplicity, we drop the $^{`}+^{`}$
sign in the superscript.
\par
Choose a larger stencil $T=\{I_{i-r},...,I_{i+r}\}$. Consider  a fourth degree
polynomial($r=2$)	 based on the nodal point information of the numerical flux which satisfies
\begin{equation}
\frac{1}{\Delta x}\int_{I_{j}}p(\xi)d\xi=\frac{1}{\Delta x}\int_{I_{j}}h(\xi)d\xi=f_{j}, j=i-r,...,i+r.
\end{equation}
Evaluating   this polynomial $p(x)$ at $x=x_{i+\frac{1}{2}}$ gives
\begin{equation}
\hat{f}_{i+\frac{1}{2}}:=P(x_{i+\frac{1}{2}})  = \frac{1}{60}(2f_{i-2}-13f_{i-1}+47f_{i}+27f_{i+1}-3f_{i+2}).
\end{equation}
\par
If there is a discontinuity inside the stencil $T$, then the corresponding interpolation process to the approximation of flux $\hat{f}_{i+\frac{1}{2}}$ may generate oscillations. In order to alleviate this the WENO procedure is employed, in which the stencil $T$ is divided into $(r+1)$ smaller stencils: $S_{k}=\{I_{i-r+k},...,I_{i+k}\},k=0,...,r$. The second degree polynomials $p^{k}(x),k=0,...,r$ are constructed in the associated stencils $S_{k}$ to approximate the function $h(x)$ that satisfies
\begin{equation}\nonumber
\frac{1}{\Delta x}\int_{I_{j}}p^{k}(\xi)d\xi=\frac{1}{\Delta x}\int_{I_{j}}h(\xi)d\xi=f_{j}, j=i-r+k,...,i+k; k=0,...,r.
\end{equation} 
The explicit expressions of polynomials $p^{k}(x),k=0,1,2$ as
\begin{equation*}\nonumber
\begin{aligned}
p^{0}(x)=&\frac{1}{24}\left[(-f_{i-2}+2f_{i-1}+23f_{i})+12(f_{i-2}-4f_{i-1}+3f_{i})\eta+12(f_{i-2}-2f_{i-1}+f_{i}) \eta^{2}\right],\\
p^{1}(x)=&\frac{1}{24}\left[(-f_{i-1}+26f_{i-1}-f_{i+1})+12(f_{i+1}-f_{i-1}) \eta +12(f_{i-1}-2f_{i}+f_{i+1})\eta^{2}\right],\\
p^{2}(x)=&\frac{1}{24}\left[(23f_{i}+2f_{i+1}-f_{i+2})+12(-3f_{i}+4f_{i+1}-f_{i+2})\eta+12(f_{i}-2f_{i+1}+f_{i+2})\eta^{2}\right],
\end{aligned}
\end{equation*}
where $\eta=\left(\frac{x-x_{i}}{\Delta x}\right)$. 
The evaluation of these polynomials $p^{k}(x),k=0,1,2$  at $x=x_{i+\frac{1}{2}}$ gives
\begin{eqnarray}
\begin{aligned}\label{flux5eq}
\hat{f}_{i+\frac{1}{2}}^{0}  = & \frac{1}{6}(2f_{i-2}-7f_{i-1}+11f_{i}),\\
\hat{f}_{i+\frac{1}{2}}^{1}  = & \frac{1}{6}(-f_{i-1}+5f_{i}+2f_{i+1}),\\
\hat{f}_{i+\frac{1}{2}}^{2}  = & \frac{1}{6}(2f_{i}+5f_{i+1}-f_{i+2}).
\end{aligned}
\end{eqnarray}
The Taylor's expansion of (\ref{flux5eq}) reveals
\begin{eqnarray*}\label{TF}
\hat{f}_{i+\frac{1}{2}}^{0} & = & h_{i+\frac{1}{2}}-\frac{\Delta x^{3}}{4}f{}^{(3)}(0)+O\left(\Delta x^{4}\right),\\
\hat{f}_{i+\frac{1}{2}}^{1} & = & h_{i+\frac{1}{2}}+\frac{\Delta x^{3}}{12}f{}^{(3)}(0)+O\left(\Delta x^{4}\right),\\
\hat{f}_{i+\frac{1}{2}}^{2} & = & h_{i+\frac{1}{2}}-\frac{\Delta x^{3}}{12}f{}^{(3)}(0)+O\left(\Delta x^{4}\right).
\end{eqnarray*}
The values of the function $p(x)$ at the point $x=x_{i+\frac{1}{2}}$ of cell $I_{i}$, can be written as  a linear combination of $p^{k}(x)$ at the point $x=x_{i+\frac{1}{2}}$ in the smooth regions. Thus the linear/ideal weights are defined as
\begin{equation}\label{s1}
\hat{f}_{i+\frac{1}{2}}=\sum_{k=0}^{2}d_k \hat{f}^{k}_{i+ \frac{1}{2}}.
\end{equation}
The values of these linear weights are  $d_{0}=\frac{1}{10}$, $d_{1}=\frac{3}{5}$, $d_{2}=\frac{3}{10}$. Note that each $d_k\geq0$ and $\displaystyle \sum_{k=0}^{2}d_k=1.$
\par
In the non-smooth regions,  (\ref{s1}) is not valid to  approximate the  flux function $\hat{f}_{i+\frac{1}{2}}$ in terms of local information. This issue is resolved by introducing the nonlinear weights $\omega_k$ such that
\begin{equation}\label{ew10}
\hat{f}_{i+\frac{1}{2}}=\sum_{k=0}^{2}\omega_k \hat{f}^{k}_{i+ \frac{1}{2}}.
\end{equation}
\par
 These nonlinear weights constructed in subsequent steps are such that in  smooth regions, the nonlinear weights should converge to the linear weights with the required order of accuracy and in the non-smooth regions, these have to tend to zero so that the contribution from the non-smooth regions to the approximation of the flux $\hat{f}_{i+\frac{1}{2}}$ is negligible, with this the final reconstruction is essentially non-oscillatory. Thus, the  nonlinear weights have to satisfy the following  properties: \\
{ \emph{Convexity:}}
\begin{equation}
\sum_{k=0}^{2}\omega_{k}=1,\omega_{k}\geq 0,k=0,1,2.
\end{equation}
{\emph{ Optimal Order:}}
 If $f$ is smooth
in  stencil $T$, then
\begin{equation}
L(u)=f^{\prime}(x_{i})+O\bigg(\Delta x^{5}\bigg).
\end{equation}
{ \emph{ENO property:}}
If a substencil $T^{D}\subset T$
contains a discontinuity of $f$, but there exists another sub-stencil
$T^{C}\subset T$ where $f$ is smooth,
then
\begin{eqnarray*}
\begin{aligned}
\omega_{D} & =O(\Delta x^{q})\text{ {for\,some} }q>0, \text{and}\\
\omega_{C} & =\varTheta(1),
\end{aligned}
\end{eqnarray*}
as $\Delta x\rightarrow 0,$ where $O(\cdot)$ and $\varTheta(\cdot)$ are standard Bachmann-Landau notation \cite{DB}. 
\par
The following result 
relate the effective order of accuracy of a WENO scheme to the difference
between its non-linear weights $\omega_{k}$ and the linear weights $d_{k}.$
\begin{lem}
(\textbf{Sufficient Condition}) If the nonlinear weights satisfy the condition
\begin{eqnarray}
\begin{cases}
\omega_{k}-d_{k} & =O(\Delta x^{2}),\\
\omega_{k}^{+}-\omega_{k}^{-} & =O(\Delta x^{3}),
\end{cases}k=0,1,2,
\end{eqnarray}
or
\begin{equation}\label{20}
\omega_{k}-d_{k}=O(\Delta x^{3}),k=0,1,2,
\end{equation}
then the corresponding WENO scheme satisfy the optimal order of accuracy, where the superscripts $'+'$ or $'-'$ on $\omega_{k}$ correspond to their use in  $\hat{f}_{i+ \frac{1}{2}}$ or $\hat{f}_{i- \frac{1}{2}}$ respectively \cite{henrick aslam powers5, Borges caramona10}.
\end{lem}
\subsection{WENO-LOC weights and its order of convergence}
The nonlinear weights  defined in \cite{XD Liu} are 
\begin{equation}\label{WLOC1}
\omega_k=\frac{\alpha_k}{\displaystyle\sum_{l=0}^{2}\alpha_l}, \, \alpha_k=\frac{d_k}{\left(\epsilon+\beta_{k}\right)^{p}},
\end{equation}
where  $\epsilon$ is a small positive number which is set to be $\epsilon=10^{-5}$ to avoid division by zero, $p=2$ is chosen to increase the difference of scales of distinct weights at non-smooth parts of the solution. Note that $\alpha_k$ are the unnormalized weights and $\omega_k$ are the normalized weights. The  smoothness of the flux is measured by the derivatives of the reconstructed flux $\hat{f}_{i+\frac{1}{2}}^{k}$ on each stencil $S_k$, $k=0,1,2,$ based on the undivided differences as 
\begin{eqnarray}
\beta_{k}=\sum_{n=1}^{2}\sum_{m=1}^{3-n}\frac{(f[i+k+m-3,n])^{2}}{3-n}, k=0,1,2,
\end{eqnarray}
where $f[\cdot,\cdot]$ is the $n^{\text{th}}$ undivided difference,
\begin{eqnarray*}
f[i,0]&=&f_i, \\
f[i,n]&=&f[i+1,n-1]-f[i,n-1].\nonumber
\end{eqnarray*}
So, we have
\begin{equation}
\beta_k = \frac{1}{2}\left( (f{[i+k-2,1]})^2+(f{[i+k-1,1]})^2)\right)+(f{[i+k-2,2]})^2,k=0,1,2,
\end{equation}
and its explicit form for $k=0,1,2$ are
\begin{eqnarray}\label{r3UDsmooth}
\begin{aligned}
\beta_0 &= \frac{1}{2}\left( (f_{i-1}-f_{i-2})^2+(f_{i}-f_{i-1})^2) \right)+(f_{i}-2f_{i-1}+f_{i-2})^2,\\
\beta_1 &= \frac{1}{2}\left((f_{i}-f_{i-1})^2+(f_{i+1}-f_{i})^2)\right)+(f_{i-1}-2f_{i}+f_{i+1})^2,\\
\beta_2 &= \frac{1}{2}\left((f_{i+1}-f_{i})^2+(f_{i+2}-f_{i+1})^2)\right)+(f_{i+2}-2f_{i+1}+f_{i})^2.
\end{aligned}
\end{eqnarray}
The Taylor's  expansion of the smoothness indicator \eqref{r3UDsmooth} of the candidate stencils  at $x=x_i$ are expressed as
\begin{align}
\beta_{0} & =\left(f_{i}^{\prime}\right)^{2}\Delta x^{2}-2f_{i}^{\prime}f_{i}^{\prime \prime}\Delta x^{3}+\left(\frac{4}{3}f_{i}^{\prime}f_{i}^{\prime\prime\prime}+\frac{9}{4}\left(f_{i}^{\prime\prime}\right)^{2}\right)\Delta x^{4}-\frac{23}{6}f_{i}^{\prime\prime}f_{i}^{\prime\prime\prime}\Delta x^{5}+O(\Delta x^{6}),\nonumber \\
\beta_{1} & =\left(f_{i}^{\prime}\right)^{2}\Delta x^{2}+\left(\frac{1}{3}f_{i}^{\prime}f_{i}^{\prime\prime\prime}+\frac{5}{4}\left(f_{i}^{\prime\prime}\right)^{2}\right)\Delta x^{4}+O(\Delta x^{6}), \label{eq:Tbetak}\\
\beta_{2} & =\left(f_{i}^{\prime}\right)^{2}\Delta x^{2}+2f_{i}^{\prime}f_{i}^{\prime\prime}\Delta x^{3}+\left(\frac{4}{3}f_{i}^{\prime}f_{i}^{\prime\prime\prime}+\frac{9}{4}\left(f_{i}^{\prime\prime}\right)^{2}\right)\Delta x^{4}+\frac{23}{6}f_{i}^{\prime\prime}f_{i}^{\prime\prime\prime}\Delta x^{5}+O(\Delta x^{6}),\nonumber 
\end{align}
Substituting  \eqref{eq:Tbetak} into \eqref{WLOC1}, we get
\begin{eqnarray*}
\alpha_0&=&\frac{1}{10\left(f_{i}^{\prime}\right)^{4}\Delta x^{4}}\left(1+4\frac{f_{i}^{\prime\prime}}{f_{i}^{\prime}}\Delta x-\left(\frac{8}{3}\frac{f_{i}^{\prime\prime\prime}}{f_{i}^{\prime}}+\frac{9}{2}\frac{\left(f_{i}^{\prime\prime}\right)^{2}}{\left(f_{i}^{\prime}\right)^{2}}\right)\Delta x^{2}+O(\Delta x^3)\right),\\
\alpha_1&=&\frac{6}{10\left(f_{i}^{\prime}\right)^{4}\Delta x^{4}}\left(1-\left(\frac{2}{3}\frac{f_{i}^{\prime\prime\prime}}{f_{i}^{\prime}}+\frac{5}{2}\frac{\left(f_{i}^{\prime\prime}\right)^{2}}{\left(f_{i}^{\prime}\right)^{2}}\right)\Delta x^{2}+O(\Delta x^3)\right),\\
\alpha_2&=&\frac{3}{10\left(f_{i}^{\prime}\right)^{4}\Delta x^{4}}\left(1-4\frac{f_{i}^{\prime\prime}}{f_{i}^{\prime}}\Delta x-\left(\frac{8}{3}\frac{f_{i}^{\prime\prime\prime}}{f_{i}^{\prime}}+\frac{9}{2}\frac{\left(f_{i}^{\prime\prime}\right)^{2}}{\left(f_{i}^{\prime}\right)^{2}}\right)\Delta x^{2}+O(\Delta x^3)\right),
\end{eqnarray*}
and
\begin{eqnarray}
\omega_0&=& \frac{1}{10}+\frac{2f_{i}^{\prime\prime}}{5f_{i}^{\prime}}\Delta x+O(\Delta x^2),\nonumber \\
\omega_1&=& \frac{6}{10}+O(\Delta x^4),\label{WL1} \\ 
\omega_2&=& \frac{3}{10}-\frac{6f_{i}^{\prime\prime}}{5f_{i}^{\prime}}\Delta x+O(\Delta x^2).\nonumber
\end{eqnarray}
\par Note that in the above procedure a small parameter $\epsilon$ is omitted since it is only used to avoid the denominator to be zero. From  \eqref{WL1}, it is concluded that the nonlinear weights approaches to the linear weights with  first order of accuracy. So, the numerical scheme with WENO-LOC weights provides the overall fourth order of accuracy  in smooth regions and further the order of accuracy degrades to third-order in presence of first-order critical points(which can observed by doing similar analysis on the unnormalized and normalized weights).   
\subsection{WENO-JS weights and its order of convergence}
\par
As observed in above, the smoothness indicators of WENO-LOC scheme does not achieve the optimal order of convergence in the smooth regions, the authors Jiang and Shu in \cite{Jiang and shu7} constructed a new smoothness measurements $\beta_k$ in (\ref{WLOC1}) based on the concept of reducing the total variation of the numerical solution on  each stencil as,
\begin{equation}\label{21}
\beta_{k}=\sum_{l=1}^{2}\Delta x^{2l-1}\intop_{x_{i-\frac{1}{2}}}^{x_{i+\frac{1}{2}}}\left(\frac{d^{l}\hat{f}^{k}}{dx^{q}}\right)^{2}dx,
\end{equation}
which is a scaled square sum of all the derivatives of  interpolation polynomial $\hat{f}^{k}(x)$ over the interval
 $\left(x_{i-\frac{1}{2}},x_{i+\frac{1}{2}}\right)$. The explicit form of these smoothness indicators are as follows
\begin{eqnarray}
\beta_{0} & = & \frac{13}{12}(f_{i-2}-2f_{i-1}+f_{i})^{2}+\frac{1}{4}(f_{i-2}-4f_{i-1}+3f_{i})^{2},\nonumber \\
\beta_{1} & = & \frac{13}{12}(f_{i-1}-2f_{i}+f_{i+1})^{2}+\frac{1}{4}(f_{i+1}-f_{i-1})^{2},\label{eq:22}\\
\beta_{2} & = & \frac{13}{12}(f_{i}-2f_{i+1}+f_{i+2})^{2}+\frac{1}{4}(3f_{i}-4f_{i+1}+f_{i+2})^{2}.\nonumber
\end{eqnarray}
By Taylor's expansion of these smoothness indicators, one can
obtain
\begin{align}
\beta_{0} & =\left(f_{i}^{\prime}\right)^{2}\Delta x^{2}+\left(\frac{13}{12}\left(f_{i}^{\prime\prime}\right)^{2}-\frac{2}{3}f_{i}^{\prime}f_{i}^{\prime\prime\prime}\right)\Delta x^{4}+\left(\frac{-13}{6}f_{i}^{\prime\prime}f_{i}^{\prime\prime\prime}+\frac{1}{2}f_{i}^{\prime}f_{i}^{iv}\right)\Delta x^{5}+O(\Delta x^{6}),\nonumber \\
\beta_{1} & =\left(f_{i}^{\prime}\right)^{2}\Delta x^{2}+\left(\frac{13}{12}\left(f_{i}^{\prime\prime}\right)^{2}+\frac{1}{3}f_{i}^{\prime}f_{i}^{\prime\prime\prime}\right)\Delta x^{4}+O(\Delta x^{6}),\label{eq:tb}\\
\beta_{2} & =\left(f_{i}^{\prime}\right)^{2}\Delta x^{2}+\left(\frac{13}{12}\left(f_{i}^{\prime\prime}\right)^{2}-\frac{2}{3}f_{i}^{\prime}f_{i}^{\prime\prime\prime}\right)\Delta x^{4}+\left(\frac{13}{6}f_{i}^{\prime\prime}f_{}^{\prime\prime\prime}-\frac{1}{2}f_{i}^{\prime}f_{i}^{iv}\right)\Delta x^{5}+O(\Delta x^{6}).\nonumber
\end{align} 
\par 
 Now, let us see the order of convergence of the nonlinear weights of WENO-JS5 scheme. Substituting  \eqref{eq:tb} into \eqref{WLOC1} with $p=2$ and $\epsilon=0$, we get
\begin{eqnarray*}
\alpha_0&=&\frac{1}{\left(f_{i}^{\prime}\right)^{4}\Delta x^{4}}\left(\frac{1}{10}-\left(\frac{-2}{15}\frac{f_{i}^{\prime\prime\prime}}{f_{i}^{\prime}}+\frac{13}{60}\frac{\left(f_{i}^{\prime\prime}\right)^{2}}{\left(f_{i}^{\prime}\right)^{2}}\right)\Delta x^{2}+O(\Delta x^3)\right),\\
\alpha_1&=&\frac{1}{\left(f_{i}^{\prime}\right)^{4}\Delta x^{4}}\left(\frac{6}{10}-\left(\frac{2}{5}\frac{f_{i}^{\prime\prime\prime}}{f_{i}^{\prime}}+\frac{13}{10}\frac{\left(f_{i}^{\prime\prime}\right)^{2}}{\left(f_{i}^{\prime}\right)^{2}}\right)\Delta x^{2}+O(\Delta x^3)\right),\\
\alpha_2&=&\frac{1}{\left(f_{i}^{\prime}\right)^{4}\Delta x^{4}}\left(\frac{3}{10}-\left(\frac{-2}{5}\frac{f_{i}^{\prime\prime\prime}}{f_{i}^{\prime}}+\frac{13}{20}\frac{\left(f_{i}^{\prime\prime}\right)^{2}}{\left(f_{i}^{\prime}\right)^{2}}\right)\Delta x^{2}+O(\Delta x^3)\right),
\end{eqnarray*}
and
\begin{eqnarray}
\omega_0&=& \frac{1}{10}+\frac{3f_{i}^{\prime\prime}}{25f_{i}^{\prime}}\Delta x^2+O(\Delta x^3),\nonumber \\
\omega_1&=& \frac{6}{10}-\frac{12f_{i}^{\prime\prime}}{25f_{i}^{\prime}}\Delta x^2+O(\Delta x^3),\label{WJS12} \\ 
\omega_2&=& \frac{3}{10}+\frac{9f_{i}^{\prime\prime}}{25f_{i}^{\prime}}\Delta x^2+O(\Delta x^3).\nonumber
\end{eqnarray}
\par
From \eqref{WJS12}, we conclude that the WENO-JS5 nonlinear weights converges to the ideal weights with the second-order of accuracy. So, the numerical scheme with WENO-JS5 weights provides the overall fifth order of accuracy  in smooth regions. Note that the advantage with the WENO-JS5 weights over the WENO-LOC weights is that it improves the one order of accuracy in smooth regions.  Further the order of accuracy of WENO-JS5 scheme degrades to third-order in presence of first-order critical points and to second-order if the second derivatives vanishes.
\section{Construction of a new nonlinear weights}
A novel global smoothness measurement is constructed based on the linear combination of undivided differences of second-order derivatives which leads to provide a sixth-order of accuracy on the global stencil $S^{5}$ as 
\begin{eqnarray}\label{r3g}\nonumber
\begin{aligned}
\zeta=\bigg|\left(\left(f_{i-2}-2f_{i-1}+f_{i}\right)^2-2\left(f_{i-1}-2f_{i}+f_{i+1}\right)^2+\left(f_{i}-2f_{i+1}+f_{i+2}\right)^2\right)\bigg|,
\end{aligned}
\end{eqnarray}
and the Taylor's expansion of $\zeta$ gives
\begin{equation}\label{r3gt}\nonumber
\zeta = 2\bigg|\left(f_i{''}f_i^{(4)}+(f_i{'''})^2\right)\bigg|\Delta x^6+O(\Delta x^8).
\end{equation}
Note that in the construction of  WENO-LOC weights, the  usage of the first-order derivatives in the smoothness indicators are not able to produce the required order of accuracy i.e., third-order,  because of this reason, we avoid the first-order derivatives information in the construction of global smoothness measurement of the global stencil.
Now, we define the nonlinear weights $\omega_k$ as
\begin{equation}\label{weightsud5}
\omega_k=\frac{\alpha_k}{\displaystyle\sum_{k=0}^{2}\alpha_k},k=0,1,2,
\end{equation}
and the unnormalized weights as
\begin{equation}\label{alphaud}
\alpha_k=d_k\left(1+\dfrac{\zeta}{\beta_k+\epsilon}\right),\;k=0,1,2,
\end{equation}
such that the nonlinear weights $\omega_k$ converge to the ideal weights with the higher order of accuracy where we use \eqref{r3UDsmooth} the smoothness indicators   $\beta_k$. The parameter $\epsilon$ is taken as a small number  to avoid the division by zero and chosen this value as $10^{-16}$. 
\par
Now, we check the convergence order of nonlinear weights in smooth regions i.e., $f_{i}^{\prime}\neq0$. Substituting  \eqref{eq:Tbetak}  into \eqref{alphaud}, we have
\begin{eqnarray*}
\alpha_0&=&\frac{1}{10}+\frac{1}{5}\left(\frac{f_{i}^{\prime\prime}f_{i}^{iv}}{\left(f_{i}^{\prime}\right)^{2}}+\frac{\left(f_{i}^{\prime\prime\prime}\right)^{2}}{\left(f_{i}^{\prime}\right)^{2}}\right)\Delta x^{4}+O(\Delta x^5),\\
\alpha_1&=&\frac{6}{10}+\frac{6}{5}\left(\frac{f_{i}^{\prime\prime}f_{i}^{iv}}{\left(f_{i}^{\prime}\right)^{2}}+\frac{\left(f_{i}^{\prime\prime\prime}\right)^{2}}{\left(f_{i}^{\prime}\right)^{2}}\right)\Delta x^{4}+O(\Delta x^5),\\
\alpha_2&=&\frac{3}{10}+\frac{3}{5}\left(\frac{f_{i}^{\prime\prime}f_{i}^{iv}}{\left(f_{i}^{\prime}\right)^{2}}+\frac{\left(f_{i}^{\prime\prime\prime}\right)^{2}}{\left(f_{i}^{\prime}\right)^{2}}\right)\Delta x^{4}+O(\Delta x^5),
\end{eqnarray*}
\begin{eqnarray}
\sum_{k=0}^{2}\alpha_{k}=\left(1+2\left(\frac{f_{i}^{\prime\prime}f_{i}^{iv}}{\left(f_{i}^{\prime}\right)^{2}}+\frac{\left(f_{i}^{\prime\prime\prime}\right)^{2}}{\left(f_{i}^{\prime}\right)^{2}}\right)\Delta x^{4}+O(\Delta x^5)\right),
\end{eqnarray}
and now from \eqref{weightsud5}
\begin{eqnarray}
\omega_0= \frac{1}{10}+O(\Delta x^5), \,
\omega_1= \frac{6}{10}+O(\Delta x^5),\label{WU1} \,
\omega_2= \frac{3}{10}+O(\Delta x^5).
\end{eqnarray}
From \eqref{WU1}, the proposed nonlinear weights converges to the ideal weights with the fifth-order of accuracy \eqref{20}. 
\par
 Now, we analyze the nonlinear weights  \eqref{weightsud5} with \eqref{alphaud} in presence of first-order critical points. From  \eqref{eq:Tbetak},   \eqref{weightsud5} and \eqref{alphaud} with $f^{\prime}_{i}=0$, we have,
\begin{eqnarray}
\omega_0&=& \frac{1}{10}+\frac{752}{1125}\left(\frac{f_{i}^{iv}}{f_{i}^{\prime\prime}}+\frac{\left(f_{i}^{\prime\prime\prime}\right)^{2}}{\left(f_{i}^{\prime\prime}\right)^{2}}\right)\Delta x^{2}+O(\Delta x^3),\nonumber \\
\omega_1&=& \frac{6}{10}+\frac{64}{375}\left(\frac{f_{i}^{iv}}{f_{i}^{\prime\prime}}+\frac{\left(f_{i}^{\prime\prime\prime}\right)^{2}}{\left(f_{i}^{\prime\prime}\right)^{2}}\right)\Delta x^{2}+O(\Delta x^3),\label{WUD12} \\ 
\omega_2&=& \frac{3}{10}-\frac{16}{125}\left(\frac{f_{i}^{iv}}{f_{i}^{\prime\prime}}+\frac{\left(f_{i}^{\prime\prime\prime}\right)^{2}}{\left(f_{i}^{\prime\prime}\right)^{2}}\right)\Delta x^{2}+O(\Delta x^3).\nonumber
\end{eqnarray}
So, at the first-order critical points, the sufficient condition \eqref{20} is not satisfied, as it resembles the numerical scheme can not achieve the desired fifth-order accuracy. To achieve the desired order of accuracy in presence of first-order critical points, we define the  unnormalized weights by introducing a parameter $p$ as
\begin{equation}\label{alphaud1}
\alpha_k=d_k\left(1+\left(\dfrac{\zeta}{\beta_k+\epsilon}\right)^{p}\right),\;k=0,1,2,
\end{equation}
and considered this $p$ value as $2$ which is an integer, so the nonlinear weights achieves the desired sufficient condition \eqref{20}.  To confirm this, we analyze the weight by substituting  \eqref{eq:Tbetak}  into \eqref{weightsud5} with \eqref{alphaud1}, we get
\begin{eqnarray}
\omega_0&=& \frac{1}{10}-\frac{681}{686}\left(\frac{f_{i}^{iv}}{f_{i}^{\prime\prime}}+\frac{\left(f_{i}^{\prime\prime\prime}\right)^{2}}{\left(f_{i}^{\prime\prime}\right)^{2}}\right)^{2}\Delta x^{4}+O(\Delta x^5),\nonumber \\
\omega_1&=& \frac{6}{10}+\frac{761}{2369}\left(\frac{f_{i}^{iv}}{f_{i}^{\prime\prime}}+\frac{\left(f_{i}^{\prime\prime\prime}\right)^{2}}{\left(f_{i}^{\prime\prime}\right)^{2}}\right)^{2}\Delta x^{4}+O(\Delta x^5), \\ 
\omega_2&=& \frac{3}{10}-\frac{457}{1995}\left(\frac{f_{i}^{iv}}{f_{i}^{\prime\prime}}+\frac{\left(f_{i}^{\prime\prime\prime}\right)^{2}}{\left(f_{i}^{\prime\prime}\right)^{2}}\right)^{2}\Delta x^{4}+O(\Delta x^5).\nonumber
\end{eqnarray}
\par
At the first-order critical points with $p=1$, the nonlinear weights do not achieves the desired fifth-order accuracy. Note that the nonlinear weights $\omega_{k},k=0,1,2$ satisfies the convexity property and also achieving the optimal order in smooth regions. Now, we check the ENO-property for the proposed nonlinear weights.
\subsection{ENO-property for proposed nonlinear weights}
 The proposed nonlinear weights are 
\begin{equation}\label{WLOC}
\omega_k=\frac{\alpha_k}{\displaystyle\sum_{l=0}^{2}\alpha_l}, \, \alpha_k=d_k\left(1+\left(\dfrac{\zeta}{\beta_k+\epsilon}\right)^{p}\right).
\end{equation}
If a sub-stencil  $T^{C}\subset T$ which is smooth then the smoothness indicators and unnormalized weights are 
\begin{eqnarray*}
\begin{aligned}
\beta_{C} & =O(\Delta x^{2}), \,
\alpha_{C} & =d_k\left(1+\left(\dfrac{\zeta}{\beta_C+\epsilon}\right)^{p}\right),
\end{aligned}
\end{eqnarray*} 
 and if a sub-stecnil  $T^{D}\subset T$ is discontinuous then 
\begin{eqnarray*}
\begin{aligned}
\beta_{D} & =\varTheta(1), \,
\alpha_{D} & =d_k\left(1+\left(\dfrac{\zeta}{\beta_D+\epsilon}\right)^{p}\right),
\end{aligned}
\end{eqnarray*} 
where $\epsilon$ is not  predominant factor and note that $\zeta=\varTheta(1)$. 
Now, 
\begin{eqnarray}\nonumber
\frac{\alpha_D}{\alpha_C}&=&\displaystyle\frac{d_D\left(1+\left(\dfrac{\zeta}{\beta_D+\epsilon}\right)^{p}\right)}{d_C\left(1+\left(\dfrac{\zeta}{\beta_C+\epsilon}\right)^{p}\right)},\\ \nonumber
&=& \displaystyle\frac{d_D\left({\beta_C+\epsilon}\right){^p}\left(\zeta^p+\left({\beta_D+\epsilon}\right){^p}\right)}{d_C\left({\beta_D+\epsilon}\right){^p}\left(\zeta^p+\left({\beta_c+\epsilon}\right){^p}\right)},\\ \nonumber
\end{eqnarray}
\begin{eqnarray}\nonumber
&=&\displaystyle\frac{d_D\left(\left(\dfrac{\beta_C+\epsilon}{\zeta}\right)^{p}+\left(\dfrac{\beta_C+\epsilon}{\beta_D+\epsilon}\right)^{p}\right)}{d_C\left(1+\left(\dfrac{\beta_C+\epsilon}{\zeta}\right)^{p}\right)},\\ \nonumber
&=& \varTheta(1) \displaystyle \frac{\left(\left(\dfrac{O(\Delta x^{2})+\epsilon}{\varTheta(1)}\right)^{p}+\left(\dfrac{O(\Delta x^{2})+\epsilon}{\varTheta(1)}\right)^{p}\right)}{\left(1+\left(\dfrac{O(\Delta x^{2})+\epsilon}{\varTheta(1)}\right)^{p}\right)},\\ \nonumber
&=&O(\Delta x^{2p}). \nonumber
\end{eqnarray}
So, 
\begin{eqnarray}\nonumber
\omega_D=\begin{cases}
O(\Delta x^{2}), & \text{if} \, \, p=1,\\
O(\Delta x^{4}), & \text{if} \, \, p=2,
\end{cases}
\hspace{0.5cm} \omega_C=\varTheta(1), 
\end{eqnarray}
which concludes that as the mesh is refining the weight assigned to the discontinuous stencil  $\omega_D$ tends to zero, so the defined nonlinear weights satisfies the ENO-property. Note that for $p=1$, the weight assigned to the discontinuous stencil is larger in comparison to the  weight assigned to the discontinuous stencil for $p=2$ case. Now, we conduct a  test which contains the smooth regions and critical points to check the same numerically. First, we check the convergence order in smooth regions and later we check for the critical point case. 
\subsection{Accuracy test}
\textbf{Case 1:}
Consider the linear advection equation,
\begin{equation}\label{lint}
u_t+u_x=0,\;-1\leq x \leq1 ,t>0,
\end{equation}
with  a smooth initial data
\begin{equation}\label{sin11}
  u_0(x)=\sin(\pi x).
\end{equation}
\par
We employed periodic boundary conditions and evaluated up to time $t=2$ to verify the order of convergence.  The $L_{1}$ and $L_{\infty}$-errors along with their numerical order of convergence is calculated with the WENO-LOC, WENO-JS5 and the proposed scheme (hereafter we call it as WENO-UD5 scheme). Note that, we use fourth-order non TVD Runge-Kutta method by the time step $\Delta t\approx \Delta x^{5/4}$ which is effectively fifth-order. The value of  $p$ is considered as  $2$ for WENO-LOC and WENO-JS5 schemes whereas for the WENO-UD5 scheme, we use $p=1$ and $p=2$ to verify the order of convergence.
\begin{table}[ht!]\label{table2}
	\footnotesize
\centering
   \begin{tabular}{c *{9}{c}}
\toprule
    N &
    \multicolumn{2}{p{3cm}}{WENO-LOC} &
    \multicolumn{2}{p{3cm}}{WENO-JS5} &
    \multicolumn{2}{p{4cm}}{WENO-UD5(p=1)} &
    \multicolumn{2}{p{4cm}}{WENO-UD5(p=2)} & \\
 \cmidrule(lr){2-3} \cmidrule(lr){4-5} \cmidrule(lr){6-7} \cmidrule(lr){8-9} 
    &
    $L_1$-error& $L_1$-order&
     $L_1$-error&$ L_1$-order&
      $L_1$-error&$ L_1$-order&
       $L_1$-error&$ L_1$-order&\\
 \midrule
    10  &8.9483e-03&    -& 3.0143e-02&-& 5.3749e-03&	-& 6.2259e-03&-\\
    20  &1.8809e-03&2.2502&1.4794e-03&4.3487&2.0589e-04&4.7063&2.1028e-04&	4.8879\\
    40 &2.8548e-04&	2.7200&4.5012e-05&5.0386&6.5442e-06&	4.9755&6.5629e-06&	5.0018\\
    80 & 2.0902e-05&3.7717& 1.3984e-06&5.0085&2.0340e-07&	5.0078&2.0345e-07&	5.0116\\
    160 & 1.2973e-06&4.0101&4.3604e-08&5.0032&6.3301e-09&5.0059&6.3302e-09&	5.0063\\
    320 & 5.3161e-08&4.6090&1.3598e-09&5.0030&1.9741e-10&5.0030&1.9741e-10&	5.0030 \\
    640 & 1.0097e-09&5.7184& 4.2207e-11&5.0098&6.1851e-12&4.9963&6.1851e-12&	4.9963\\
     \bottomrule
    &
    $L_\infty$-error& $L_\infty$-order&
     $L_\infty$-error&$ L_\infty$-order&
     $L_\infty$-error&$ L_\infty$-order&
     $L_\infty$-error&$ L_\infty$-order&\\
 \cmidrule(lr){2-3} \cmidrule(lr){4-5} \cmidrule(lr){6-7} \cmidrule(lr){8-9}
    10  &1.2594e-02&    -&  4.8506e-02&	-&8.1305e-03&	-&1.0439e-02&\\
    20  &4.3976e-03&1.5179&2.5414e-03&4.2545&3.5455e-04&	4.5193&3.3755e-04&4.9507 \\
    40  &7.2682e-04&2.5970&8.9204e-05&4.8324&1.1745e-05&4.9159&1.0291e-05&5.0356 \\
    80 & 8.2692e-05&3.1358& 2.7766e-06&5.0057&3.6572e-07&5.0052&3.1904e-07&5.0115\\
    160 &8.7588e-06&3.2389&8.6040e-08&5.0122&1.0971e-08&5.0590&9.9414e-09&5.0041\\
    320 &6.1263e-07&3.8376&2.5528e-09&5.0749&3.2988e-10&5.0556& 3.1008e-10&5.0027\\
    640 &1.2653e-08&5.5975&7.3502e-11&5.1182&1.0066e-11&5.0344&9.7160e-12&4.9961\\
 \bottomrule
\end{tabular}
   \caption{ $L_1$ and $L_\infty$-error and orders with initial condition \eqref{sin11}.}
   \label{Tab1}
\end{table}
From the numerical errors and its order of convergence from table \eqref{Tab1}, it concludes that the WENO-LOC scheme converges to fourth-order of accuracy. Note that, the WENO-LOC scheme commits the lesser error on fine mesh  which results a super-convergence phenomena. As per the case of WENO-JS5 scheme, it  converges to fifth-order accuracy  and when it comes to WENO-UD5 schemes it achieves the fifth-order of accuracy. The advantage of WENO-UD5 scheme over the WENO-JS5 scheme  is that the WENO-UD5 scheme produce very lesser errors especially on the coarser mesh.\\  
\textbf{Case 2:}
In this case, we use the  initial condition
\begin{equation}\label{sin21}
u_0(x)=\sin\left(\pi x-\frac{\sin(\pi x)}{\pi}\right),
\end{equation}
for the equation \eqref{lint} which contains first-order critical point i.e., $u_x =0$ in [-1,1] but  $u_{xxx} \neq 0$. This test case has its own importance in the literature  since it has been shown that WENO-JS5 scheme do not achieve the desired order of convergence rate at critical point case and  as a result it attracted to many of researchers of this field.   We calculated the  numerical errors and its order of convergence for the WENO-LOC, WENO-JS5 schemes with the parameter $p=2$,  for the WENO-UD5 scheme with $p=1,2$ and are tabulated in \eqref{Tab2}.
\begin{table}[ht!]
	\footnotesize
\centering
   \begin{tabular}{c *{9}{c}}
\toprule
    N &
    \multicolumn{2}{p{3cm}}{WENO-LOC} &
    \multicolumn{2}{p{3cm}}{WENO-JS5} &    
    \multicolumn{2}{p{4cm}}{WENO-UD5(p=1)} &
    \multicolumn{2}{p{4cm}}{WENO-UD5(p=2)} &\\
 \cmidrule(lr){2-3} \cmidrule(lr){4-5} \cmidrule(lr){6-7} \cmidrule(lr){8-9}
    &
    $L_1$-error& $L_1$-order&
    $L_1$-error&$ L_1$-order&
    $L_1$-error&$ L_1$-order&
     $L_1$-error&$ L_1$-order&\\
 \midrule
    10  &5.6764e-02&	-&6.1696e-02&-&6.3213e-02&	-&4.0544e-02&-\\		
    20  &8.0131e-03&2.8245&4.9323e-03&	3.6448&2.6393e-03&4.5820&2.0967e-03&	4.2733\\
    40 &1.5787e-03&	2.3436&3.6462e-04&	3.7578&7.8995e-05&	5.0623&7.4596e-05&	4.8129\\
    80 &2.2341e-04&	2.8210&1.7098e-05&	4.4145&2.4010e-06&	5.0401&2.3500e-06&	4.9884\\
    160 &2.4940e-05&3.1632&7.3414e-07&	4.5416&7.4563e-08&5.0090&7.3372e-08&	5.0013\\
    320 &1.3371e-06&4.2213&2.5134e-08&	4.8683&2.3268e-09&5.0020&2.2907e-09&	5.0014\\
    640 &1.5086e-08&6.4698&5.1179e-10&	5.6179&7.2648e-11&5.0013&7.1514e-11&	5.0014\\
 \bottomrule
    &
    $L_\infty$-error& $L_\infty$-order&
     $L_\infty$error&$ L_\infty$-order&
      $L_\infty$error&$ L_\infty$-order&
     $L_\infty$error&$ L_\infty$-order&\\
 \cmidrule(lr){2-3} \cmidrule(lr){4-5} \cmidrule(lr){6-7} \cmidrule(lr){8-9}
    10 &1.1502e-01&	-&1.3639e-01&-&1.3294e-01&	-&8.1286e-02&-\\
    20 &2.0379e-02&	2.4967&1.2790e-02&3.4146&6.9116e-03&	4.2656&5.0463e-03&4.0097\\
    40 &5.0171e-03&2.0222&1.0952e-03&3.5458&2.2836e-04&4.9196&2.1071e-04&4.5819\\
    80 &9.7368e-04&2.3653&8.7557e-05&3.6448&6.6880e-06&5.0936&6.7014e-06&4.9747\\
    160 &1.6311e-04&2.5776&7.4148e-06&3.5617&2.0989e-07&4.9939&2.0988e-07&4.9968\\
    320 &1.7206e-05&3.2449&4.0271e-07&4.2026&6.5526e-09&5.0014&6.5526e-09&5.0014\\
   640 &2.9112e-07&5.8852&6.4373e-09&5.9671&2.0485e-10&4.9994&2.0485e-10&4.9994\\
 \bottomrule
\end{tabular}
   \caption{$L_1$ and $L_\infty$-error and orders with initial condition (\ref{sin21}).}
   \label{Tab2}
\end{table}
\par
It is shown that the increase in the parameter value $p$, the theoretical order of convergence  achieves for the WENO-UD5 scheme and as a result  the numerical scheme gets desired fifth-order of convergence. 
\par
Note that from the table \eqref{Tab2} at first-order critical points, the proposed WENO-UD5 scheme achieves the desired fifth-order of convergence for the parameter value $p=1$ too. As an  immediate consequence that  it is giving an intuition about the nonlinear weights such that these are not necessarily have to satisfy the sufficient condition. But whether it really reflects in achieving the required ENO order with the designed weights for the parameter $p=1$. To know this, now we analyze the ENO property in the numerics for the WENO-LOC, WENO-JS5,  WENO-UD5 schemes.   That is, does  WENO-UD5 reconstruction satisfies the ENO order in presence of discontinuities? or does the numerical scheme with the designed nonlinear weights achieves atleast the ENO-order of accuracy as WENO-LOC or  WENO-JS5 scheme in presence of discontinuities?
\subsection{Reconstruction in the  discontinuous case}
To analyze the  nonlinear weights in presence of discontinuities, let us define 
\begin{equation}
K=\{k:f\text{ is not smooth in }S^{k}\textbf{ }\}
\end{equation}
then  for WENO-LOC and WENO-JS5 schemes, we have
\begin{eqnarray}
\omega_{k}=\frac{\alpha_{k}}{\displaystyle\sum_{l=0}^{2}\alpha_{l}}=\begin{cases}
O\big(\Delta x^{2}\big) & \text{ if }\,\,k\in K,\\
\varTheta(1) & \text{ if }\,\,k\not\in K.
\end{cases}
\end{eqnarray}
Therefore,
\begin{eqnarray}
\begin{aligned}
f\Big(x_{i+\frac{1}{2}}\Big)-\hat{f}\Big(x_{i+\frac{1}{2}}\Big) & =f\Big(x_{i+\frac{1}{2}}\Big)-\sum_{k=0}^{2}\omega_{k}\hat{f}_{i+\frac{1}{2}}^{k},\\
 & =\sum_{k=0}^{2}\omega_{k}\bigg(f\Big(x_{i+\frac{1}{2}}\Big)-\hat{f}_{i+\frac{1}{2}}^{k}\bigg),\\
 & =\sum_{k\notin K}\omega_{k}\bigg(f\Big(x_{i+\frac{1}{2}}\Big)-\hat{f}_{i+\frac{1}{2}}^{k}\bigg)+\sum_{k\in K}\omega_{k}\bigg(f\Big(x_{i+\frac{1}{2}}\Big)-\hat{f}_{i+\frac{1}{2}}^{k}\bigg),\\
 & =\sum_{k\notin K}\varTheta(1)O(\Delta x^{3})+\sum_{k\in K}O(\Delta x^{2})\varTheta(1),\\
 & =O(\Delta x^{2}).
\end{aligned}
\end{eqnarray}
Thus,  the order of accuracy of WENO-LOC and WENO-JS5 scheme is worse than the corresponding ENO scheme which is of the order $3.$ Now, we show the ENO-order for proposed WENO-UD5 scheme. For this, we have
\begin{eqnarray}
\omega_{k}=\frac{\alpha_{k}}{\displaystyle\sum_{l=0}^{2}\alpha_{l}}=\begin{cases}
O\big(\Delta x^{2p}\big) & \text{ if }\,\,k\in K,\\
\varTheta(1) & \text{ if }\,\,k\not\in K.
\end{cases}
\end{eqnarray}
Thus,
\begin{eqnarray*}
\begin{aligned}
f\Big(x_{i+\frac{1}{2}}\Big)-\hat{f}\Big(x_{i+\frac{1}{2}}\Big) & =f\Big(x_{i+\frac{1}{2}}\Big)-\sum_{k=0}^{2}\omega_{k}\hat{f}_{i+\frac{1}{2}}^{k},\\
 & =\sum_{k=0}^{2}\omega_{k}\bigg(f\Big(x_{i+\frac{1}{2}}\Big)-\hat{f}_{i+\frac{1}{2}}^{k}\bigg),\\
 & =\sum_{k\notin K}\omega_{k}\bigg(f\Big(x_{i+\frac{1}{2}}\Big)-\hat{f}_{i+\frac{1}{2}}^{k}\bigg)+\sum_{k\in K}\omega_{k}\bigg(f\Big(x_{i+\frac{1}{2}}\Big)-\hat{f}_{i+\frac{1}{2}}^{k}\bigg),\\
 & =\sum_{k\notin K}\varTheta(1)O(\Delta x^{3})+\sum_{k\in K}O(\Delta x^{2p})\varTheta(1),\\
 & =\min\left(O(\Delta x^{3}),O(\Delta x^{2p})\right).
\end{aligned}
\end{eqnarray*}
Therefore, the numerical scheme with the proposed weights have the similar behavior as WENO-LOC and WENO-JS5 schemes near the discontinuities with $p=1$ but for $p=2$, the proposed nonlinear weights achieves the desired ENO-order of accuracy i.e., atleast $3$ near the discontinuities.
\par
For more understanding of this phenomena, we analyze how the WENO reconstruction behaves in presence of discontinuities by conducting an example as follows:
consider a discontinuous function
\[
f(x)=\begin{cases}
x^{3}+\cos(x) & \text{ if }\,\,x\leq0.5,\\
x^{3}+\cos(x)+1 & \text{ if }\,\,x>0.5,
\end{cases}
\]
and a uniform grid on $[-1,1]$ with $N=\{25,50,100,200,400,800,1600\}$, note that $f^{\prime}(0)=0$.
We compute the errors of the approximations by the WENO-LOC, WENO-JS5 and WENO-UD5
reconstructions at the points $x_{i\pm1}$ where $x_{i-1}$
is at the left part of the discontinuity and $x_{i+1}$ is at the
right part of the discontinuity, $0.5\in[x_{i},x_{i+1}).$ In this
experiment, we use $\epsilon=10^{-6}$ for WENO-LOC scheme, $\epsilon=10^{-6} $ for WENO-JS5 scheme and  $\epsilon=10^{-16} $ for the WENO-UD5 schemes with the parameter $p=1$, $2$.
The results are displayed in tables \ref{Tab3}, \ref{Tab4}, \ref{Tab5} and \ref{Tab6} respectively.
We also display the deduced orders $o_{i\pm1}(\Delta x)=log_{2}(e_{i\pm1}(\Delta x/2)/e_{i\pm1}(\Delta x))$
to reveal the order of the WENO-LOC, WENO-JS5 and WENO-UD5 reconstructions.
\begin{table}[ht!]
\centering
\begin{tabular}{cccccc}
\hline
N  & $\Delta x$  & $e_{i-1}$  & $o_{i-1}$  & $e_{i+1}$  & $o_{i+1}$\tabularnewline
\hline
25  & 8.000e-02  & 6.5645e-03  & -----  & -6.2122e-03  & ------\tabularnewline
50  & 4.000e-02  & 1.8265e-03  & 1.8456  & -1.4413e-03  & 2.1077\tabularnewline
100  & 2.000e-02  & 6.4936e-04  & 1.4920  & -3.3261e-04  & 2.1155\tabularnewline
200  & 1.000e-02  & 1.7659e-04  & 1.8786  & -7.8032e-05  & 2.0917\tabularnewline
400  & 5.000e-03  & 4.5570e-05  & 1.9542 & -1.8621e-05  & 2.0671\tabularnewline
800  & 2.500e-03  & 1.1526e-05  & 1.9832  & -4.5402e-06  & 2.0361\tabularnewline
1600  & 1.250e-03  & 2.8893e-06  & 1.9963  & -1.1265e-06  & 2.0109\tabularnewline
\hline
\end{tabular}\caption{\label{Tab3} WENO-LOC with $\epsilon=10^{-6}$}
\end{table}
\begin{table}[ht!]
\centering
\begin{tabular}{cccccc}
\hline
N  & $\Delta x$  & $e_{i-1}$  & $o_{i-1}$  & $e_{i+1}$  & $o_{i+1}$\tabularnewline
\hline
25  & 8.000e-02  & 5.7174e-03  & -----  & -4.6591e-03  & ------\tabularnewline
50  & 4.000e-02  & 2.8495e-03  &1.0047  & -1.1513e-03  & 2.0168e+00\tabularnewline
100  & 2.000e-02  & 7.3615e-04  & 1.9526  & -2.8835e-04  & 1.9974e+00\tabularnewline
200  & 1.000e-02  & 1.8482e-04  & 1.9939  & -7.2085e-05  & 2.0001e+00\tabularnewline
400  & 5.000e-03  & 4.6254e-05  & 1.9985 & -1.8015e-05  & 2.0005e+00\tabularnewline
800  & 2.500e-03  & 1.1567e-05  & 1.9996  & -4.5019e-06  & 2.0006e+00\tabularnewline
1600  & 1.250e-03  & 2.8910e-06  & 2.0004  & -1.1250e-06  & 2.0006e+00\tabularnewline
\hline
\end{tabular}
\caption{\label{Tab4} WENO-JS5 with $\epsilon=10^{-6}$}
\end{table}
\begin{table}[ht!]
\centering
\begin{tabular}{cccccc}
\hline
N  & $\Delta x$  & $e_{i-1}$  & $o_{i-1}$  & $e_{i+1}$  & $o_{i+1}$\tabularnewline
\hline
25  & 8.000e-02  & 1.0620e-02  & -----  & -1.0581e-02  & ------\tabularnewline
50  & 4.000e-02  & 3.9005e-03  & 1.4451  & -2.4021e-03  & 2.1391\tabularnewline
100  & 2.000e-02  & 2.0246e-03  & 0.9460  & -6.6050e-04  & 1.7986\tabularnewline
200  & 1.000e-02  & 9.4622e-04  & 1.0974  & -2.1549e-04  & 1.6800\tabularnewline
400  & 5.000e-03  & 4.5809e-04  & 1.0465  & -7.8762e-05  & 1.4520\tabularnewline
800  & 2.500e-03  & 2.2555e-04  & 1.0222  & -3.2547e-05  & 1.2750\tabularnewline
1600  & 1.250e-03  & 1.1194e-04  & 1.0107  & -1.4616e-05  & 1.1550\tabularnewline
\hline
\end{tabular}\caption{\label{Tab5}WENO-UD5(p=1) with $\epsilon=10^{-16}$}
\end{table}
\begin{table}[ht!]
\centering
\begin{tabular}{cccccc}
\hline
N  & $\Delta x$  & $e_{i-1}$  & $o_{i-1}$  & $e_{i+1}$  & $o_{i+1}$\tabularnewline
\hline
25  & 8.000e-02  & 6.4452e-03  & -----  & -6.1681e-03  & ------\tabularnewline
50  & 4.000e-02  & 1.7801e-03  & 1.8563  & -1.4404e-03  & 2.0984\tabularnewline
100  & 2.000e-02  & 6.4325e-04  &1.4685  & -3.3316e-04  & 2.1122\tabularnewline
200  & 1.000e-02  & 1.7496e-04  & 1.8784  & -7.8546e-05  & 2.0846\tabularnewline
400  & 5.000e-03  & 4.5122e-05  & 1.9551  & -1.8885e-05  & 2.0563\tabularnewline
800  & 2.500e-03  & 1.1432e-05  & 1.9808  & -4.6150e-06  & 2.0328\tabularnewline
1600  & 1.250e-03  & 2.8757e-06  & 1.9911  & -1.1396e-06  & 2.0178\tabularnewline
\hline
\end{tabular}\caption{\label{Tab6}WENO-UD5(p=2) with $\epsilon=10^{-16}$}
\end{table}
\par
From the tables \ref{Tab3}, \ref{Tab4}, \ref{Tab5} and \ref{Tab6}, it concludes that WENO-LOC and WENO-JS5 schemes achieves second-order accuracy  whereas for WENO-UD5(p=1) scheme degrades to first-order and WENO-UD5(p=2) scheme gains the second-order of accuracy as WENO-LOC and WENO-JS5 schemes at the left  and right point of the discontinuities. 
\subsection{Consistency analysis:Optimal values of the parameters involved in nonlinear weights}
From the previous sections, the WENO-LOC, WENO-JS5 and WENO-UD5 schemes  satisfies the  ENO property if the parameter $\epsilon$ is not a predominant factor i.e., the order of the  $\epsilon$ should be as small as possible. If this value begins to dominate the smoothness indicators, what happens to the corresponding numerical scheme? And if it is the case what is the optimal order of this parameter in presence of arbitrary critical points. To know this, we have constructed a following theorem and the proof of this theorem follows as similar in article \cite{rathan2}.
\begin{thm} Let $\epsilon=\Delta x^{m}$
with $m\in\mathbb{R}^{+}$. The WENO reconstruction
of $f$ is defined by $\hat{f}(x)=\displaystyle\sum_{k=0}^{2}\omega_{k}\hat{f}^{k}(x)$
where
\[
\omega_{k}=\frac{\alpha_{k}}{\displaystyle\sum_{l=0}^{2}\alpha_{l}},\,\,\alpha_{k}=d_{k}\left(1+\left(\frac{\zeta}{{\beta_{k}}+\epsilon}\right)^{p}\right),\,\,k=0,1,2.
\]
Then with the parameters   $m\leq6-\dfrac{3}{p}$ and $p=2$, we have\\
 1. At the regions where $f$ is smooth:
\[
\begin{aligned}\omega_{k} & =d_{k}\left(1+O(\Delta x^{3})\right),\\
f\Big(x_{i+\frac{1}{2}}\Big)-\hat{f}\Big(x_{i+\frac{1}{2}}\Big) & =g\Big(x_{i+\frac{1}{2}}\Big)\Delta x^{5}+O(\Delta x^{6})
\end{aligned}
\]
for a locally Lipschitz function $g.$ 
\newline
2. If $f$ is not smooth in the stencil $T$ but
it is smooth in at least one of the sub-stencils $S_{k},k=0,1,2,$
then
\[
f\Big(x_{i+\frac{1}{2}}\Big)-\hat{f}\Big(x_{i+\frac{1}{2}}\Big)=O(\Delta x^{3}).
\]
\end{thm}
For numerical validation, we perform a numerical test  with initial profile
$u_0(x)=\sin^{3}\left(\pi x\right),$
for the equation \eqref{lint}. The initial condition  contains first and second-order critical points i.e., $u_x =0,u_{xx} =0,$ in [-1,1] but  $u_{xxx} \neq 0$. We display the $L_1-$ errors and its order of convergence for the WENO-LOC with $\epsilon=\{10^{-6},\Delta x, \Delta x^2, \Delta x^3 \}$, WENO-JS5 with  $\epsilon=\{10^{-6},\Delta x^2 \}$ and  WENO-UD5  with 
$\epsilon=\{10^{-6},10^{-16}, \Delta x,  \Delta x^2, \\ \Delta x^3, \Delta x^5 \}$.  The parameter $p=2$
 is considered in the numerical evaluation for all these schemes.
 \begin{table}[ht!]
 	\footnotesize
\centering
   \begin{tabular}{c *{9}{c}}
\toprule
    N &
    \multicolumn{2}{p{2cm}}{$\epsilon=10^{-6}$} &
    \multicolumn{2}{p{2cm}}{$\epsilon=\Delta x$} &    
    \multicolumn{2}{p{2cm}}{$\epsilon=\Delta x^2$} &
    \multicolumn{2}{p{2cm}}{$\epsilon=\Delta x^3$} &\\
 \cmidrule(lr){2-3} \cmidrule(lr){4-5} \cmidrule(lr){6-7} \cmidrule(lr){8-9}
    &
    $L_1$-error& $L_1$-order&
    $L_1$-error&$ L_1$-order&
    $L_1$-error&$ L_1$-order&
     $L_1$-error&$ L_1$-order&\\
 \midrule
    40 		 &1.6664e-02 &-          &1.8756e-03    &     -         &9.1419e-03    &	 -     &1.5998e-02  &-\\		
    80 		 &3.0586e-03 &2.4458     &7.3183e-05    &	4.6797      &5.4543e-04    &4.0670     &2.5786e-03   &	2.6332\\
    160		 &2.5932e-04 &3.5601     &1.5802e-06    &	5.5333      &2.6392e-05    &4.3692     &1.9362e-04   &	3.7353\\
    320 	 &5.5358e-06 &5.5498     &3.5929e-08    &	5.4588      &9.6170e-07    &4.7784     &9.8660e-06   &	4.2946\\
    640 	 &1.1884e-07 &5.5417     &1.0380e-09    &	5.1133      &3.1725e-08    &4.9219     &5.0890e-07   &	4.2770\\
    1280	 &1.7702e-09 &6.0690     &3.3074e-11    &	4.9720      &1.0021e-09    &4.9845     &2.5204e-08   &	4.3357\\
    2560     &2.2473e-11 &6.2996     &1.5438e-12    &	4.4211      &3.1372e-11    &4.9974     &1.1917e-09   &	4.4026\\
 \bottomrule
\end{tabular}
   \caption{WENO-LOC scheme}
   \label{Tab7}
\end{table}
\begin{table}[ht!]
\centering
   \begin{tabular}{c *{9}{c}}
\toprule
    N &
    \multicolumn{2}{p{2cm}}{$\epsilon=10^{-6}$} &   
    \multicolumn{2}{p{2cm}}{$\epsilon=\Delta x^2$} &\\
 \cmidrule(lr){2-3} \cmidrule(lr){4-5} 
    &
    $L_1$-error& $L_1$-order&
    $L_1$-error&$ L_1$-order&\\
 \midrule
    40 		 &6.0354e-03 &-          &3.3681e-03    &     -        \\		
    80 		 &9.1031e-04 &2.7290     &1.7769e-04    &	4.2445      \\
    160		 &4.8182e-05 &4.2398     &5.7473e-06    &	4.9503      \\
    320 	 &8.0849e-07 &5.8971     &1.7670e-07    &	5.0235      \\
    640 	 &1.3257e-08 &5.9304     &5.4959e-09    &	5.0068     \\
    1280	 &2.3166e-10 &5.8386     &1.7152e-10    &	5.0019      \\
    2560     &4.7348e-12 &5.6126     &5.4427e-12    &	4.9779      \\
 \bottomrule
\end{tabular}
   \caption{WENO-JS5 scheme}
   \label{Tab8}
\end{table}
\begin{table}[ht!]
\centering
   \begin{tabular}{c *{7}{c}}
\toprule
    N &
    \multicolumn{2}{p{2cm}}{$\epsilon=10^{-6}$} &
    \multicolumn{2}{p{2cm}}{$\epsilon=10^{-16}$} &
    \multicolumn{2}{p{2cm}}{$\epsilon=\Delta x$} &    \\
 \cmidrule(lr){2-3} \cmidrule(lr){4-5} \cmidrule(lr){6-7} 
    &
    $L_1$-error& $L_1$-order&
    $L_1$-error&$ L_1$-order&
    $L_1$-error&$ L_1$-order&\\
 \midrule
    40 		 &4.8292e-03 &-          &4.8412e-03    &     -         &1.1428e-03    &	 -    \\		
    80 		 &5.8219e-04 &3.0522     &6.5484e-04    &	2.8862      &3.6374e-05    &4.9735    \\
    160		 &1.7814e-06 &8.3523     &6.6937e-05    &	3.2903      &1.1389e-06    &4.9972    \\
    320 	 &3.5566e-08 &5.6464     &6.2280e-06    &	3.4260      &3.5563e-08    &5.0011    \\
    640 	 &1.1106e-09 &5.0011     &5.4925e-07    &	3.5032      &1.1106e-09    &5.0010     \\
    1280	 &3.4709e-11 &4.9999     &4.9630e-08    &	3.4682      &3.4709e-11    &4.9999     \\
    2560     &1.5177e-12 &4.5154     &5.9210e-09    &	3.0673      &1.5174e-12    &4.5156     \\
 \bottomrule
   N &
    \multicolumn{2}{p{2cm}}{$\epsilon=\Delta x^2$} &
    \multicolumn{2}{p{2cm}}{$\epsilon=\Delta x^3$} &
    \multicolumn{2}{p{2cm}}{$\epsilon=\Delta x^5$} &    \\
 \cmidrule(lr){2-3} \cmidrule(lr){4-5} \cmidrule(lr){6-7} 
    &
    $L_1$-error& $L_1$-order&
    $L_1$-error&$ L_1$-order&
    $L_1$-error&$ L_1$-order&\\
 \midrule
    40 		 &1.2097e-03 &-          &3.1734e-03    &     -         &4.8374e-03    &	 -    \\		
    80 		 &3.6389e-05 &5.0550     &7.2078e-05    &	5.4603      &6.5406e-03   &-0.4351    \\
    160		 &1.1389e-06 &4.9978     &1.2582e-06    &	5.8401      &6.6831e-05    &6.6128    \\
    320 	 &3.5563e-08 &5.0011     &3.5676e-08    &	5.1403      &6.2120e-06    &3.4274    \\
    640 	 &1.1106e-09 &5.0010     &1.1107e-09    &	5.0054      &5.4719e-07    &3.5049    \\
    1280	 &3.4709e-11 &4.9999     &3.4709e-11    &	5.0000      &4.9432e-08    &3.4685     \\
    2560     &1.5175e-12 &4.5155     &1.5176e-12    &	4.5154      &5.9363e-09    &3.0578     \\
\bottomrule
\end{tabular}
   \caption{WENO-UD5 scheme}
   \label{Tab9}
\end{table}
The tables \ref{Tab7}, \ref{Tab8} and \ref{Tab9} reveals that the WENO-LOC, WENO-JS5 schemes achieves its optimal order of accuracy in presence of critical points for the $\epsilon=\Delta x^2$ whereas  the optimal order of accuracy for the WENO-UD5 scheme achieves for the values of $\epsilon=\Delta x,\Delta x^2,\Delta x^3$. Among these, $\epsilon=\Delta x^2$ achieves globally fifth-order of accuracy with smaller errors as compared to the value of  $\epsilon=\Delta x^3$ and $\epsilon=\Delta x$. So, the conclusion for achieving the optimal order  for the parameters $\epsilon$ and $p$ in presence of arbitrary number of vanishing derivatives are $\epsilon=\Delta x^2$ and $p=2$.\\
\textbf{Computational cost:} Now we check the compuational cost of nonlinear weights for WENO-LOC and WENO-UD5 scheme.
Here $C(g)$ represents the count of $g$ and $\left\{ a\pm,b\times,c\div\right\} $
represents that number of  $a$ sums (or subtractions), $b$
products and $c$ divisions to compute $\omega_{k}.$ \\
\begin{tabular}{|c|c|c|}
		\hline 
		\textbf{Count of each parameter} & \textbf{WENO-LOC} & \textbf{WENO-UD5} \tabularnewline
		\hline 
		\hline 
		\textit{Cost per substencil} $\alpha_{k}$ & $\left\{ 1\pm,0\times,1\div\right\} +C(\beta_{k})+C(p)$ & $\left\{ 2\pm,0\times,1\div\right\} +C(\beta_{k})+C(p)+C(\zeta)$\tabularnewline
		\hline 
		\textit{Cost per stencil} $\alpha$  & $\left\{ 3\pm,0\times,3\div\right\} +3C(\beta_{k})+3C(p)$ & $\left\{ 6\pm,3\times,3\div\right\} +3C(\beta_{k})+3C(p)+C(\zeta)$\tabularnewline
		\hline 
	\textit{	Cost of $\omega$} & $C(\alpha)+\left\{ 2\pm,0\times,3\div\right\} $ & $C(\alpha)+\left\{ 2\pm,0\times,3\div\right\} $\tabularnewline
		\hline 
\end{tabular}\\

Note that $C(\alpha^{UD5})+\left\{ 2\pm,0\times,3\div\right\} =$$C(\alpha^{LOC})+\left\{ 3\pm,3\times,0\div\right\} +C(\zeta),$
so WENO-UD5 scheme is slight increase to the total cost of WENO-LOC
scheme.
\section{Numerical Results }
In this section, we have considered some benchmark problems to demonstrate the results obtained  by the proposed  scheme WENO-UD5.  For the numerical comparison purpose, we compare the results with the WENO-LOC and  WENO-JS5 schemes.  We first show the behavior of nonlinear weights by performing on a test case  i.e, we analyze how the nonlinear weights converges to the linear weights and subsequently we test the proposed scheme for the one-dimensional and two-dimensional system of Euler equations with the CFL number 0.5.
\subsection{Behaviour of nonlinear weights}
\begin{figure}[ht!]
\centering
\includegraphics[width=22cm,height=8cm]{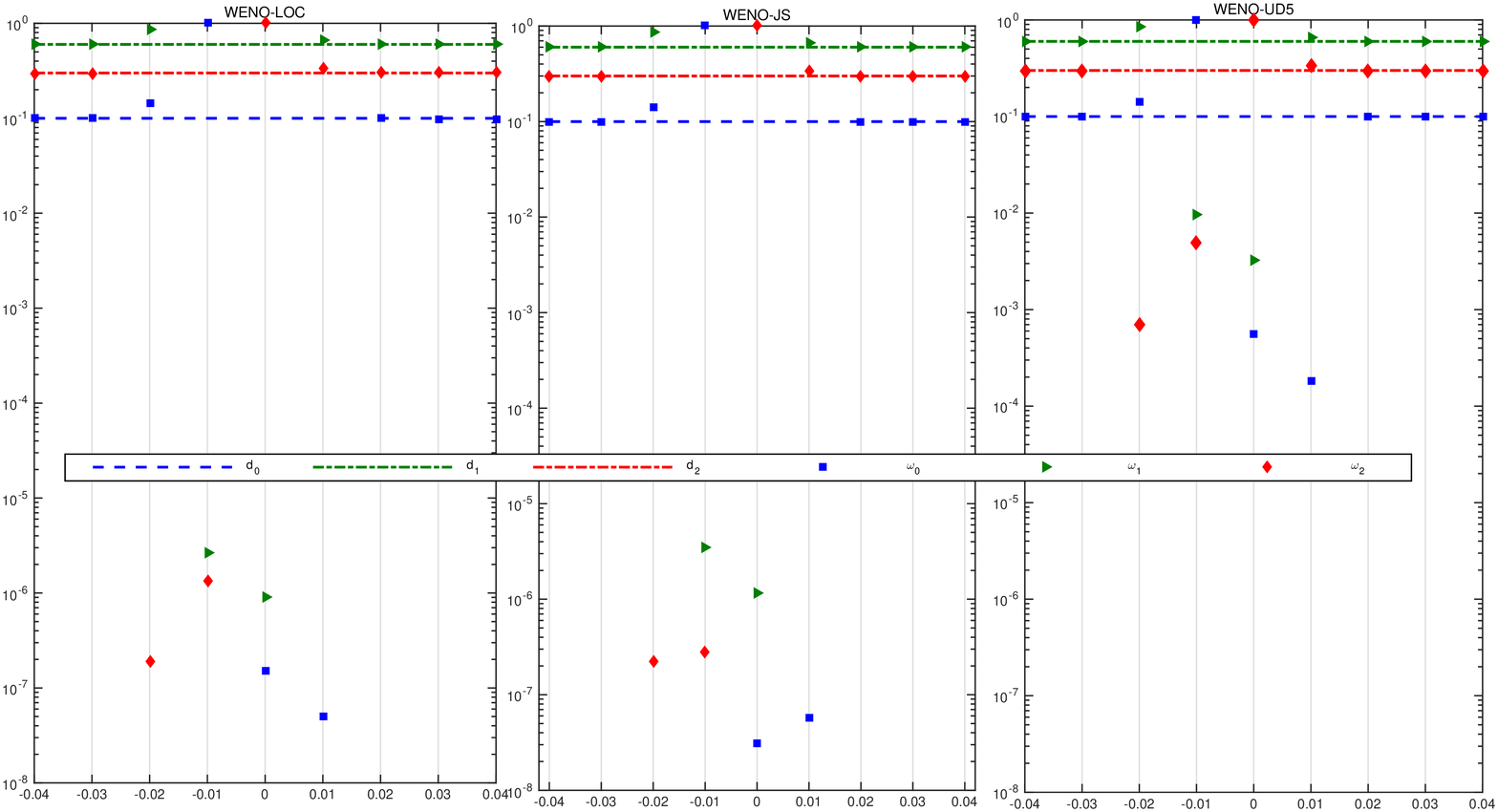}
\caption{The distribution of ideal weights $d_k$ and non-linear weights $\omega_k$, $k=0,1,2$. }\label{fig1}
\end{figure}

To understand the behavior of nonlinear weights, we considered the initial condition
\begin{eqnarray}\label{e24}
u_0(x)=
\begin{cases}
-\sin(\pi x)-\frac{1}{2}x^3, & -1<x<0,\\
-\sin(\pi x)-\frac{1}{2}x^3+1,& 0\leq x \leq 1.
\end{cases}
\end{eqnarray}
The distribution of non-linear weights $\omega_k$  and the linear weights $d_k$ are shown in Fig.\ref{fig1} for the WENO-LOC, WENO-JS5 and WENO-UD5 reconstructions. From this, it is observed that WENO-UD5 assigns larger weights for the discontinuous stencils as compared to WENO-LOC,  WENO-JS5  schemes and assigns smaller weights in smooth regions, thus the nonlinear weights are close enough to the ideal weights. 
\subsection{Linear advection test}
\begin{figure}[ht!]
\centering
\includegraphics[width=20cm,height=10cm]{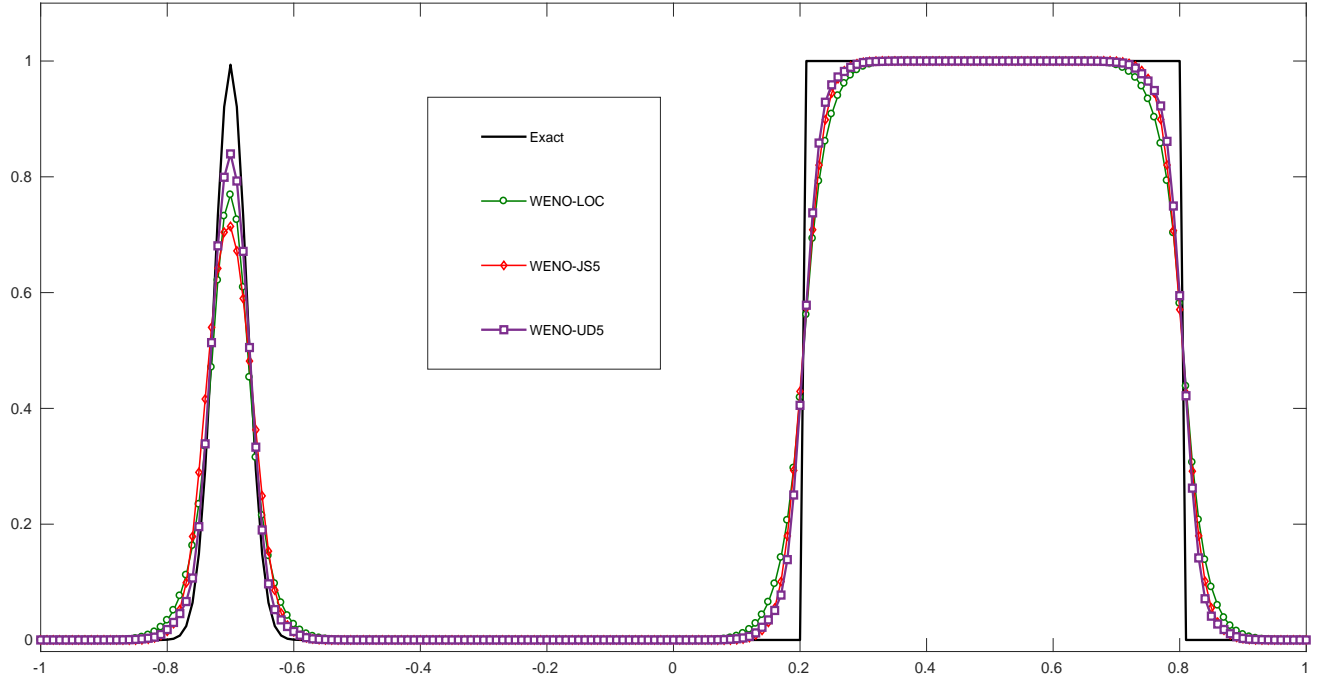}
\caption{Numerical solution of \eqref{lint} with initial condition \eqref{HRSC} }\label{fig2}
\end{figure}

Consider the linear advection equation \eqref{lint} with the  initial condition 
\begin{eqnarray}\label{HRSC}
u_0(x)=
\begin{cases}
\frac{1}{6}[G(x,z-\delta)+G(x,z)+4G(x,z+\delta)], & -0.8<x<0.2,\\
1,& 0.2\leq x \leq 0.8,
\end{cases}
\end{eqnarray}
in a computational domain $[-1,1]$ where $G(x,z)=\exp(-\beta(x-z)^2),$ $z=-0.7$ and $\beta=\frac{\log(2)}{36\delta^2}$. This initial condition consists of Gaussian and square wave shapes. The numerical solution is displayed in    Fig.\ref{fig2} for  time $t=8$. From this, it is observed  that WENO-UD5 schemes have the higher-resolution, better behavior in comparison to the WENO-LOC and WENO-JS5 schemes.
\subsection{One-dimensional Euler equations}
The numerical simulations are performed on the one-dimensional Euler equations which are given by
\begin{equation}
\begin{pmatrix}\rho\\
\rho u\\
E
\end{pmatrix}_{t}+\begin{pmatrix}\rho u\\
\rho u^{2}+p\\
u(E+p)
\end{pmatrix}_{x}=0, \label{eq:60}
\end{equation}
where $\rho,u,E,p$ are the density, velocity, total energy and pressure
respectively. The system \eqref{eq:60} represents the conservation
of mass, momentum and energy.
The total energy for an ideal polytropic
gas is defined as
\[
E=\frac{p}{\gamma-1}+\frac{1}{2}\rho u^{2},
\]
 and the eigenvalues of the Jacobian matrix $A(U)=\partial F/\partial U$ are
\begin{equation*}
  \lambda_1(u)=u-c, \quad \lambda_2(u)=u,  \quad \lambda_3(u)=u+c,
\end{equation*}
where  $U=\begin{pmatrix}\rho\\
\rho u\\
E
\end{pmatrix}$, $F=\begin{pmatrix}\rho u\\
\rho u^{2}+p\\
u(E+p)
\end{pmatrix}$ and $\gamma$ is the ratio of specific heats and its value is taken
as $1.4$.\\
\newline
Remark: For the systems of conservation laws, such as one dimensional Euler equations, the reconstruction procedures are implemented in the local characteristic directions for the purpose of avoiding spurious oscillations. For the two dimensional problems, all of these reconstruction procedures are carried out in a dimension by dimension fashion.
\subsubsection{Sod's shock  tube problem}
\begin{figure}[ht!]
\centering
\includegraphics[trim=160 10 120 10,clip,width=17cm,height=10cm]{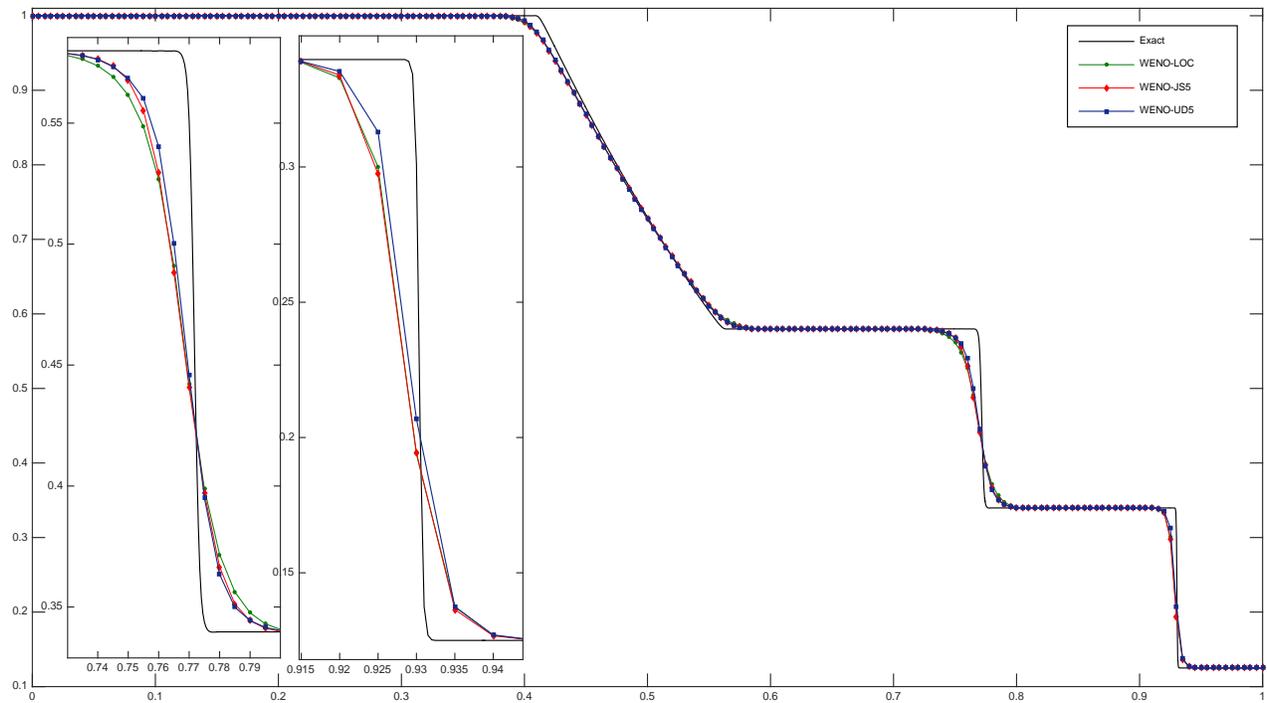}
\caption{Sod problem: Density distribution for fifth-order WENO schemes}\label{fig3}
\end{figure}

We consider the one dimensional Euler system \eqref{eq:60} with Riemann data \cite{Sod}
\[(\rho, u, p)=
\begin{cases}
(1,0,1), & -5\leq x<0,\\
(0.125,0,0.1),& 0\leq x \leq 5,
\end{cases}\]
in the computational domain $-5 \leq x \leq 5.$ The problem is initialized on the computational domain of $200$ points and is run  up to time $t=1.3,$ by this time a right-going shock wave, a right traveling contact-wave and a left-sonic  rarefaction wave establishes. The transmissive boundary conditions  are taken for numerical evaluation. The numerical results of density profiles are displayed in the  Fig.\ref{fig3}. It is observed that  the discontinuity is sharpened by  WENO-UD$5$ schemes over  WENO-LOC and WENO-JS5 schemes  due to the efficient confinement of WENO dissipation right around the discontinuity.
\subsubsection{Lax's shock tube problem}
\begin{figure}[ht!]
\centering
\includegraphics[trim=160 25 120 10,clip,width=17cm,height=10cm]{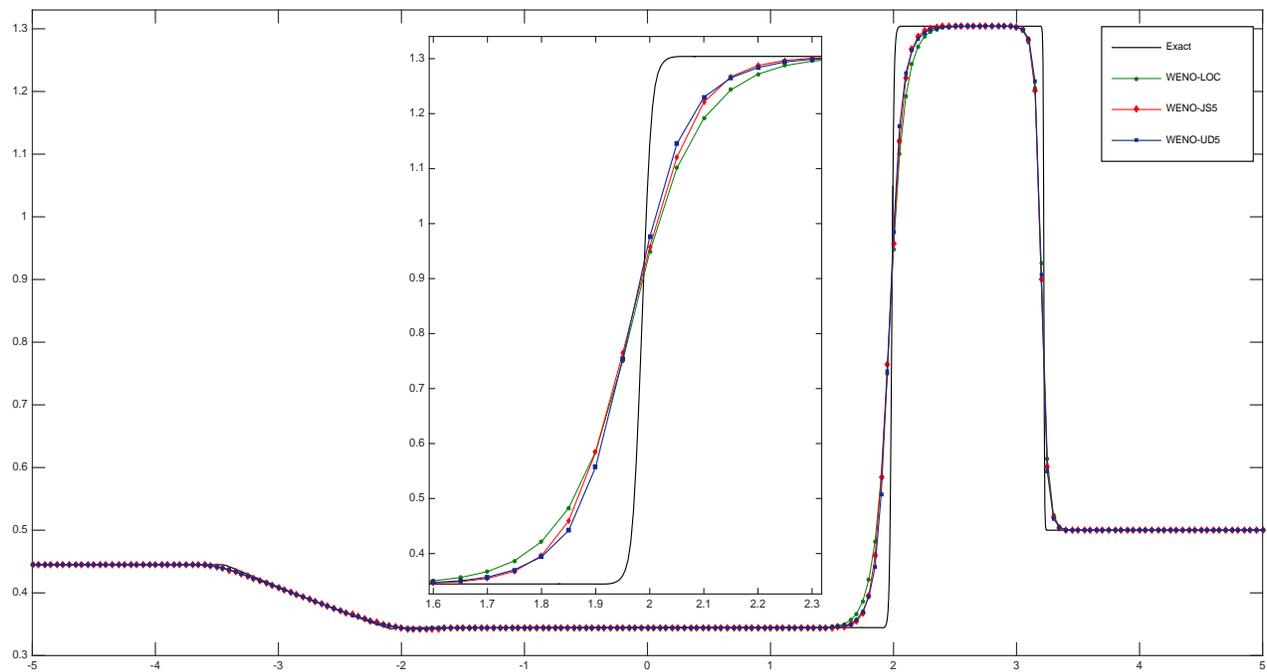}
\caption{Lax problem: Density distribution for fifth-order WENO schemes}\label{fig4}
\end{figure}

The initial condition
\[(\rho, u, p)=
\begin{cases}
(0.445,0.698,3.528), &-5\leq x<0,\\
(0.5,0,0.571),& 0\leq x \leq 5,
\end{cases}\]
is considered  \cite{PDLax} to the one dimensional Euler system of equations  \eqref{eq:60}. This shock test case is considered in the computational domain $-5 \leq x \leq 5,$  and is run up to $t=1.3$ with the zero gradient boundary conditions.  The numerical results of density profiles along with the reference solutions  are displayed in the  Fig.\ref{fig4}. The observation from the figure reveals that the numerical solutions  WENO-UD$5$ schemes have better resolution in comparison to  WENO-LOC and WENO-JS5 schemes.
\subsubsection{Mach 3 Shock entropy wave interaction test}
\begin{figure}[ht!]
\centering
\includegraphics[trim=160 30 120 10,clip,width=18cm,height=10cm]{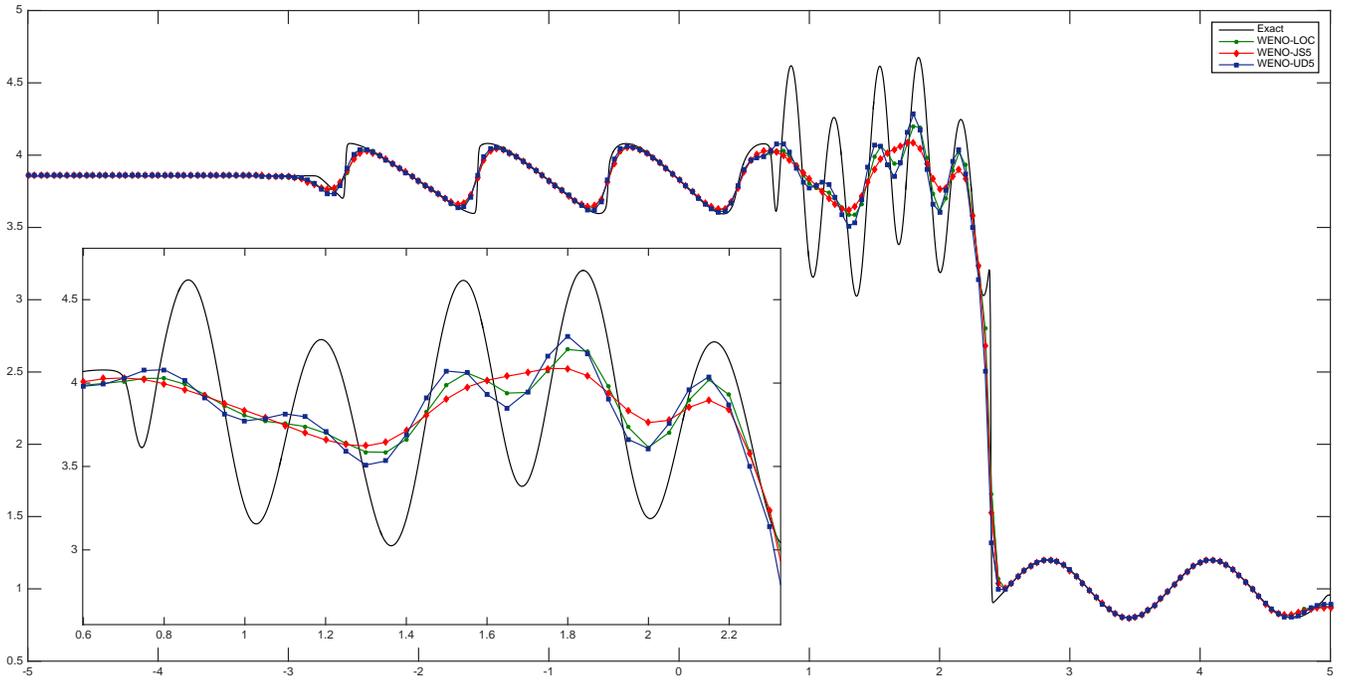}
\caption{Solution of the Mach 3 shock density wave interaction with $k=5$  with $N=200$ points.}\label{fig5}
\end{figure}

For the system \eqref{eq:60}, consider the Riemann data
\begin{equation*}
(\rho, u, p)=
  \begin{cases}
    (3.857143, 2.629369, \frac{31}{3}), & -5\leq x<-4, \\
 (1+0.2\sin(kx),0,1), & -4\leq x \leq 5,
  \end{cases}
\end{equation*}
on the spatial domain $x\in[-5,5]$ with $k=5.$   The solution of this problem  \cite{HO87} consists of a number of shocklets and fine scales structure, which are located behind a right going main shock. Fig.\ref{fig5} depicts the numerical results of WENO-LOC, WENO-JS5 and WENO-UD$5$ schemes for $N=200$ cells at time $t=1.8$ against the reference solution, computed by WENO-JS5 scheme with $N=2000$ points. We observed that WENO-UD$5$ capture more features of the solution than the WENO-LOC and WENO-JS5 particularly, at the high-frequency waves behind the shock  at deeper valleys and higher pikes in the numerical solution.
\subsection{Two dimensional Euler system of equations }
\subsubsection{2D Riemann gas dynamics problem }
\begin{figure}[ht!]
\centering
\includegraphics[trim=10 180 10 10,width=19cm,height=8cm]{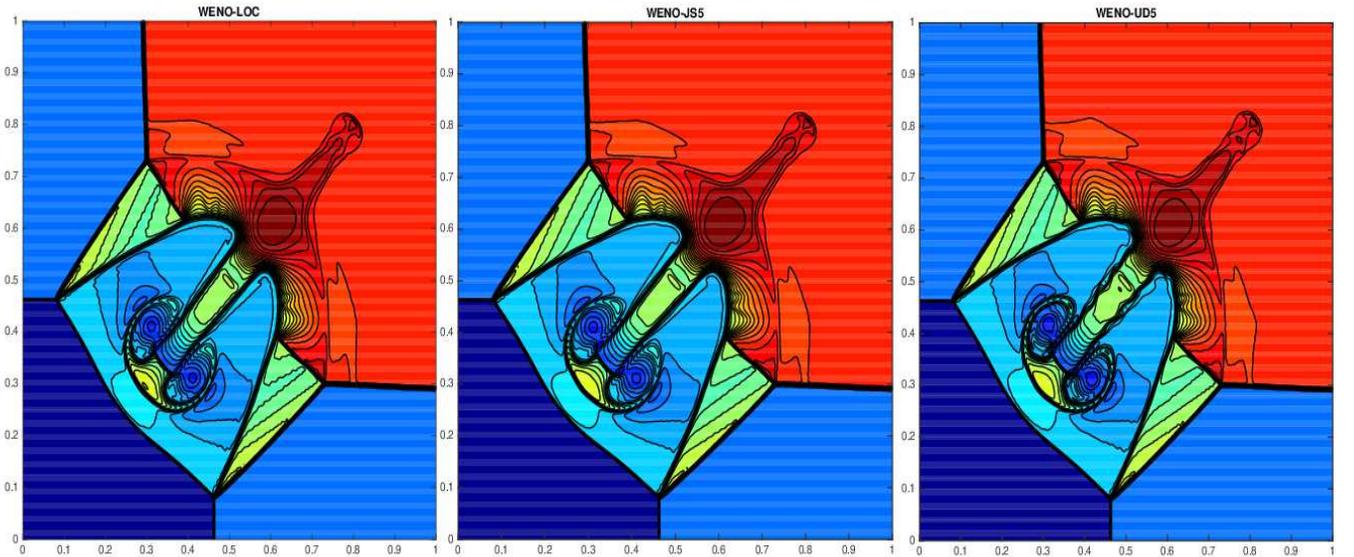}
\caption{Density profile of 2D Riemann problem of gas dynamics}\label{fig6}
\end{figure}

The two-dimensional Riemann problem of
gas dynamics \cite{Rinne} is defined by initial constant states
which is divided by the lines $x=0.8$ and $y=0.8$ on the square
as
\begin{eqnarray*}
(\rho,u,v,p) & =\begin{cases}
\begin{array}{lr}
(1.5,0,0,1.5) & \text{if}\;0.8\leq x\leq1,0.8\leq y\leq1,\\
(0.5323,1.206,0,0.3) & \text{if}\;0\leq x<0.8,0.8\leq y\leq1,\\
(0.138,1.206,1.206,0.029) & \text{if}\;0\leq x<0.8,0\leq y<0.8,\\
(0.5323,0,1.206,0.3) & \text{if}\;0.8<x\leq1,0\leq y\leq0.8,
\end{array}\end{cases}
\end{eqnarray*}
and the time evolution is governed by two-dimensional Euler equations,
\begin{eqnarray}\label{eq2d}
\begin{pmatrix}\rho\\
\rho u\\
\rho v\\
p
\end{pmatrix}_{t}+\begin{pmatrix}\rho u\\
P+\rho u^{2}\\
\rho uv\\
u(E+P)
\end{pmatrix}_{x}+\begin{pmatrix}\rho v\\
\rho uv\\
P+\rho v^{2}\\
v(E+P)
\end{pmatrix}_{y}=0.
\end{eqnarray}
The total energy $E$ and the pressure $p$ is defined by
\[
p=(\gamma-1)(E-\frac{1}{2}\rho(u^{2}+v^{2})),
\]
where $u$ and $v$ are $x$ and $y$-velocity components respectively. The
numerical solution is computed on the computational domain $[0,1]\times[0,1]$ 
with Dirichlet boundary conditions on  $400 \times 400$ grid points. According to the initial conditions,
four shocks come into being and produce a narrow jet. The numerical
solution is calculated upto time $t=0.8$.
The grid refinement results of  WENO-LOC, WENO-JS5 and WENO-UD5 schemes are given in  Fig.\ref{fig6} that contains the numerical solutions of WENO-LOC, WENO-JS5 and WENO-UD5 schemes. An examination of these results reveals that the WENO-UD$5$  scheme captures a better resolution of the fine structures  in  comparison to the schemes  WENO-LOC and WENO-JS5 respectively.
\subsubsection{Double Mach reflection of a strong shock}
\begin{figure}[ht!]
\centering
\includegraphics[width=16cm,height=5cm]{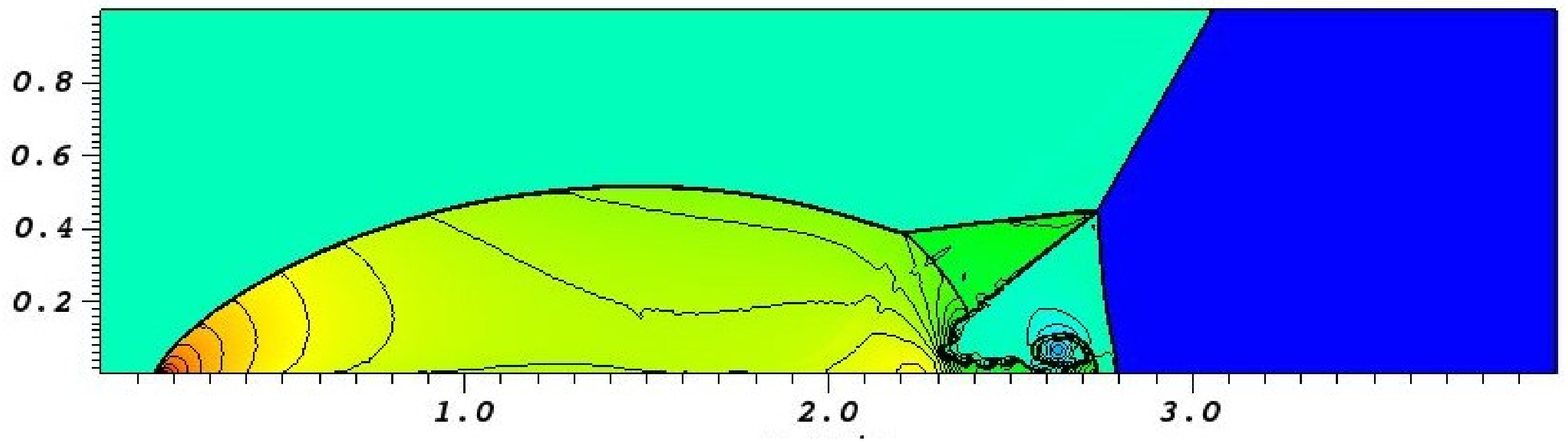}\caption{Density profile of double Mach reflection of a strong shock:WENO-LOC scheme }\label{fig18}
\includegraphics[width=16cm,height=5cm]{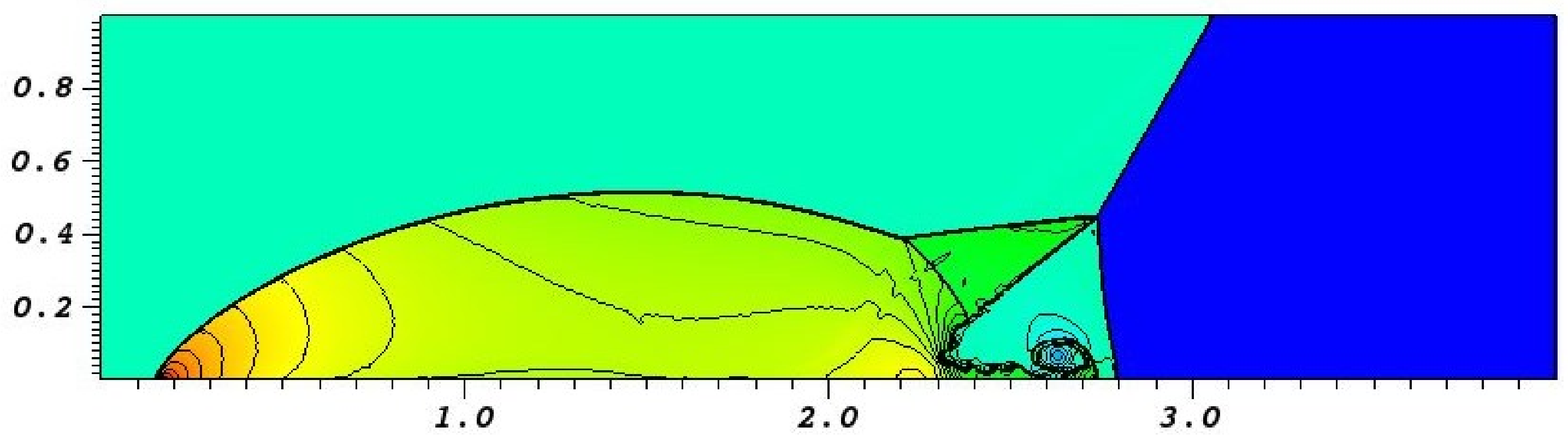}\caption{Density profile of double Mach reflection of a strong shock:WENO-JS5 scheme }\label{fig19}
\includegraphics[width=16cm,height=5cm]{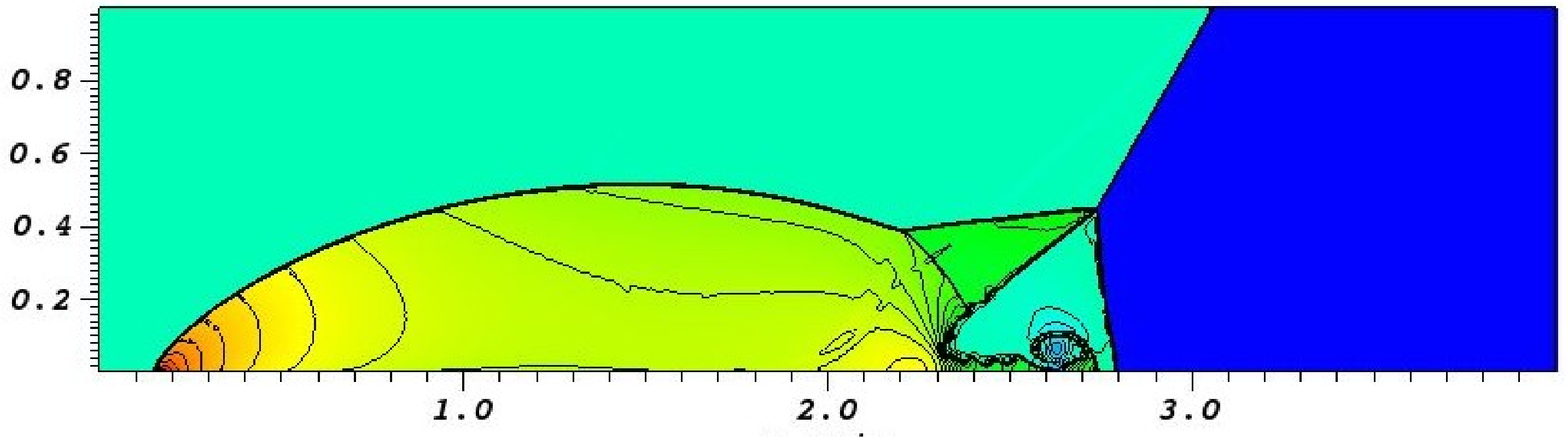}
\caption{Density profile of double Mach reflection of a strong shock:WENO-UD5 scheme }\label{fig20}
\end{figure}

\begin{figure}[ht!]
\centering
\includegraphics[width=9cm,height=7cm]{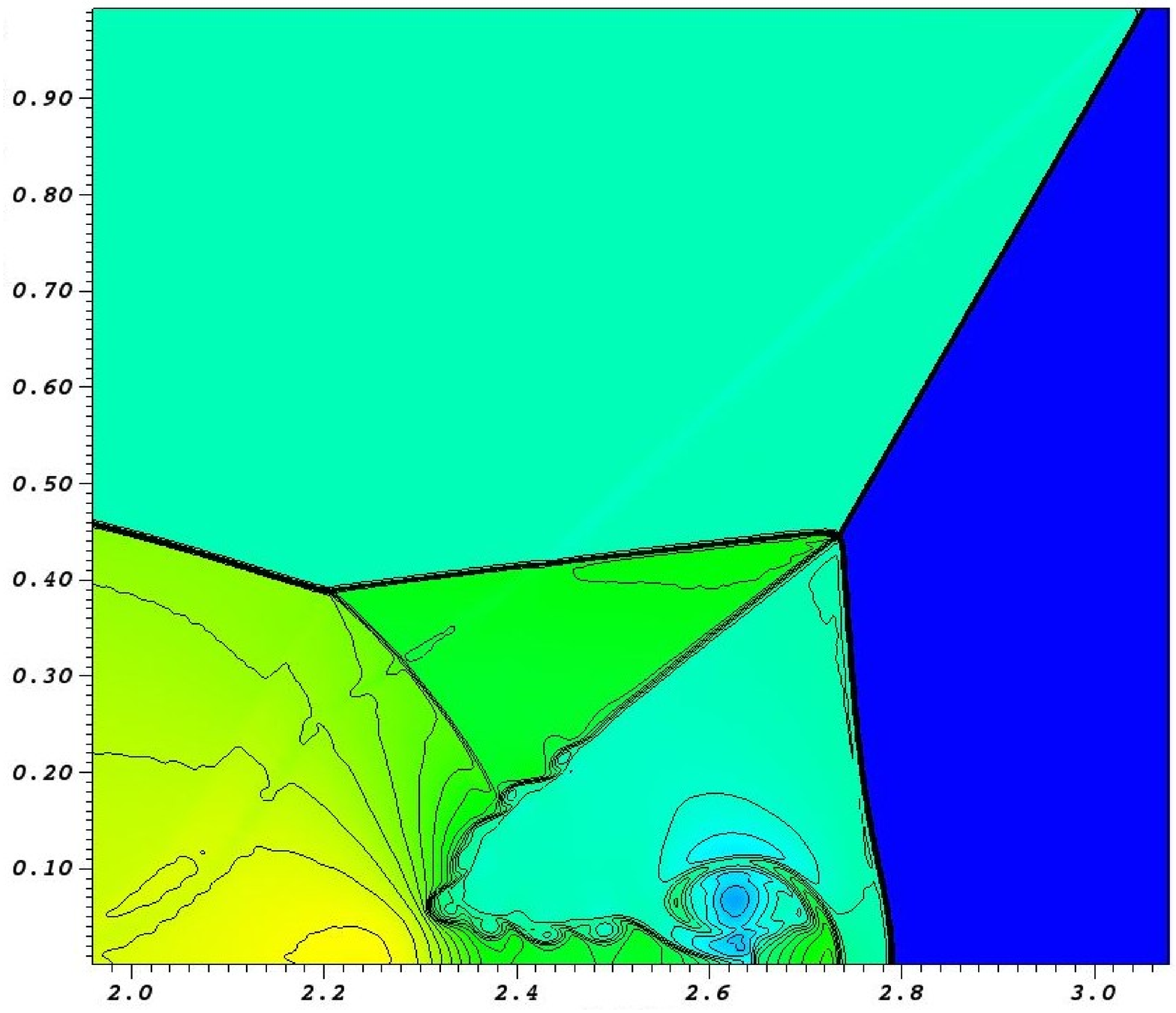}
\includegraphics[width=9cm,height=7cm]{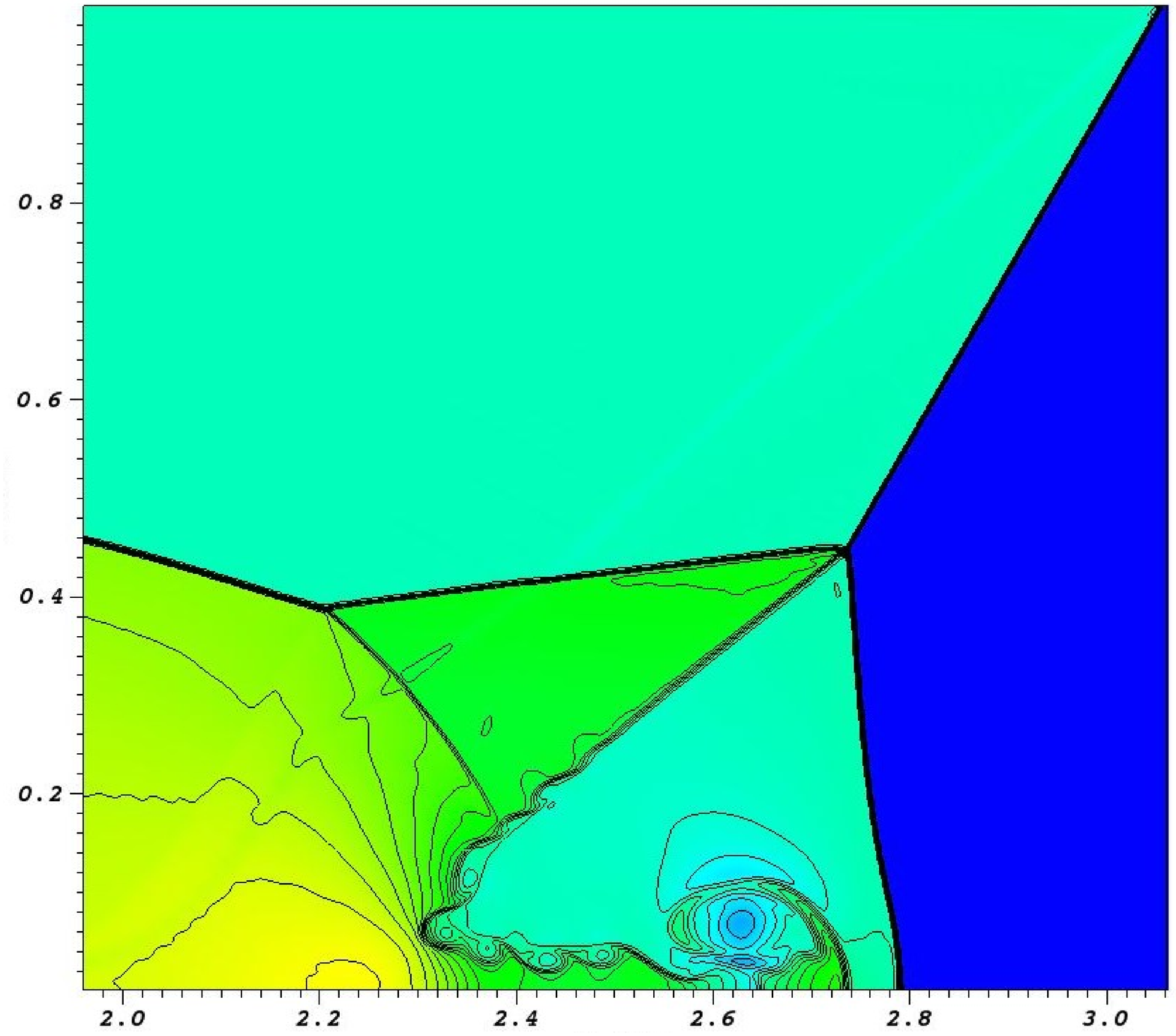}
\includegraphics[width=9cm,height=7cm]{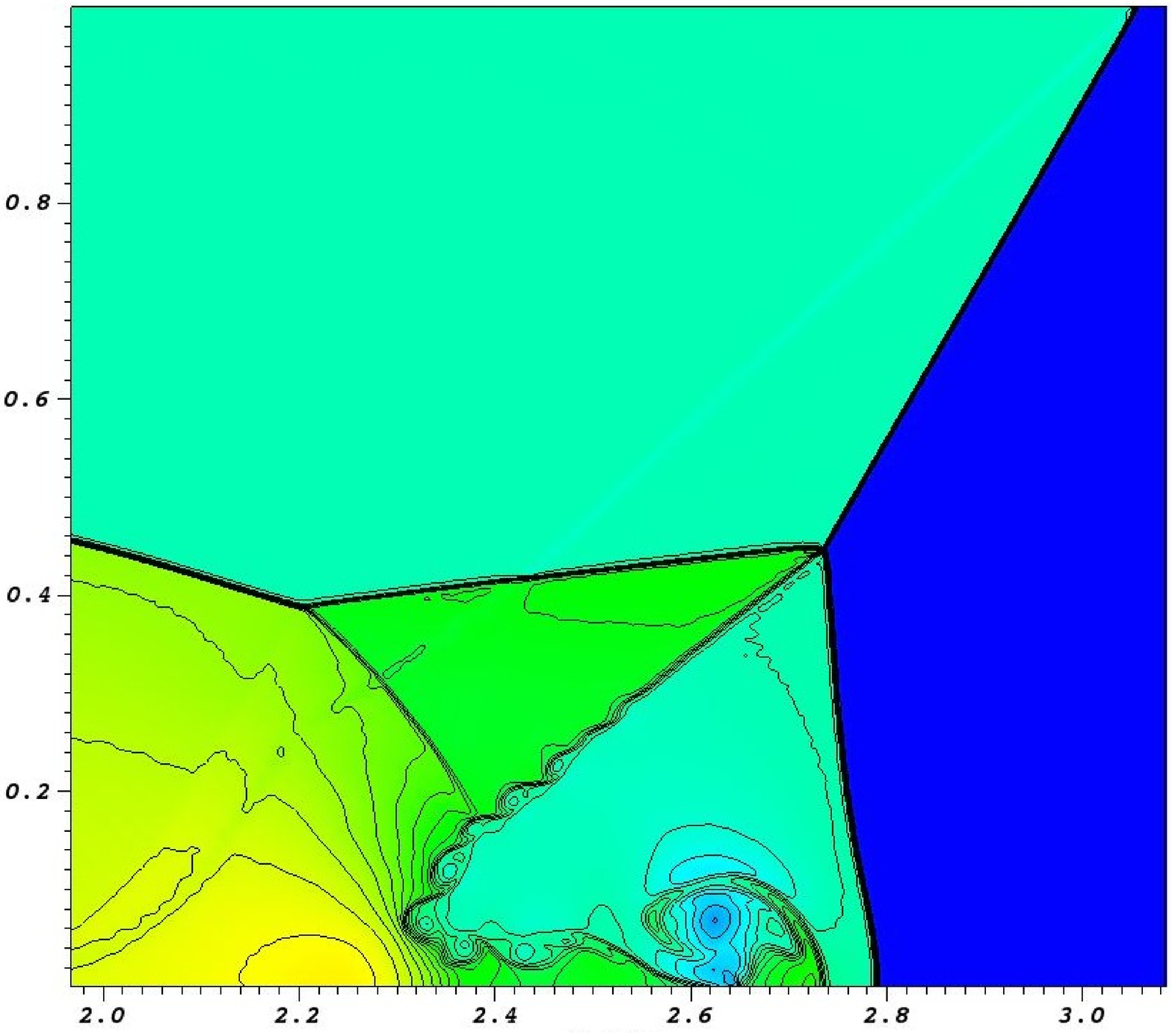}
\caption{Density profile of double Mach reflection of a strong shock around Mach stem:WENO-LOC, JS5 and UD5 scheme }\label{fig21}
\end{figure}

For the Euler's system \eqref{eq2d}, the two-dimensional
double Mach reflection problem presented by Woodward and Colella in \cite{woodward} is considered in this example on the domain of $[0,4]\times[0,1].$
The reflecting boundary conditions for $\frac{1}{6}\leq x\leq 4, y=0$ are taken. Initially, a right-moving Mach $10$ shock is positioned at $x=\frac{1}{6},y=0,$ and makes an angle of $60^{o}$
with the $x-$axis. For the bottom boundary $0\leq x< \frac{1}{6},y=0 $, the exact post shock condition is imposed. The top boundary of our computational domain uses the exact motion of the
Mach $10$ shock. Inflow and outflow boundary conditions are taken
for the left and right boundaries. The unshocked fluid has a density $1.4$,  pressure  $1.$  and  the ratio of specific heats $\gamma=1.4$. The numerical solution is computed up to time $t=0.2$ on a mesh  $500 \times 500.$ The results in the region $[0,4]\times[0,1]$ are displayed for WENO-LOC, WENO-JS5
and WENO-UD$5$ schemes are shown in Fig.\ref{fig18}, \ref{fig19} and \ref{fig20}  respectively. It can be clearly seen  in Fig.\ref{fig21} that  WENO-UD$5$ resolves the instabilities better  around the Mach stem of the problem.
\section{Conclusions}
In this paper, we have constructed a  new type of nonlinear weights for the fifth-order weighted essentially non-oscillatory scheme.  These nonlinear weights have been developed by construction of a new global smoothness indicator  using the linear combination of second-order derivative information of local stencils which resulted a sixth-order value on five point stencil. These nonlinear weights satisfies convexity, ENO property and achieves the optimal order of accuracy. The resulted numerical scheme achieved the desired fifth-order accuracy in the smooth regions and in presence of critical points. Further, we have analyzed the consistency analysis on the weight parameters  and verified the ENO property theoretically as well as numerically. Numerical results resembled  in scalar, system of one- and two-dimensional Euler equations  for typical shock tube problems and double-Mach reflection of strong shock test cases.

\end{document}